\documentclass[pdflatex,sn-mathphys-num,12pt]{sn-jnl}
\usepackage{amsmath,amsfonts,amssymb,amsthm}
\usepackage{framed}
\usepackage[dvipsnames]{xcolor}
\usepackage{graphicx}
\usepackage[T1]{fontenc}
\usepackage{times} 
\usepackage[cp1250]{inputenc}  
\usepackage{comment}
\usepackage{longtable}
\usepackage{breqn}
\usepackage{geometry}
\usepackage{mdframed}
\usepackage[T1]{fontenc} 
\definecolor{brightgreen}{rgb}{0.0, 0.5, 0.0}
\hypersetup{
    colorlinks=true,
    linkcolor=blue,
    filecolor=magenta,      
    urlcolor=blue,
    citecolor=brightgreen
} 
\makeatletter
\def\overUnderArrow{\@ifnextchar[\overUnderArrow@i{\overUnderArrow@i[]}}
\def\overUnderArrow@i[#1]#2#3{
  \ifx\relax#1\relax\array[b]{c}\overset{\text{#2}}{\uparrow}\\#3\endarray
  \else\ifx\relax#2\relax
    \array[t]{c}#3\\\underset{\text{#1}}{\downarrow}\endarray
  \else
    \array{c}\overset{\text{#2}}{\uparrow}\\#3\\\underset{\text{#1}}{\downarrow}\endarray
  \fi\fi}
\makeatother
\newcommand{\cg}{\textnormal{\textsl{g}}}
\makeatletter
\renewcommand\@makefntext[1]{%
  \noindent\makebox[10pt][r]{\@makefnmark}#1}
\makeatother

\renewcommand{\qed}{\hfill\blacksquare}
\newcommand{\del}{\partial}

\makeatother
\newtheorem{thm}{Theorem}[section]
\newtheorem{pro}[thm]{Proposition}

\theoremstyle{definition}
\newtheorem{dfn}[thm]{Definition}

\numberwithin{equation}{section}

\newcommand{\C}{\mathbb{C}} 
\renewcommand{\H}{\mathcal{H}} 
\newcommand{\T}{\mathbb{T}} 

\newbox\ncintdbox \newbox\ncinttbox
\setbox0=\hbox{$-$} \setbox2=\hbox{$\displaystyle\int$}
\setbox\ncintdbox=\hbox{\rlap{\hbox
    to \wd2{\kern-.1em\box2\relax\hfil}}\box0\kern.1em}
\setbox0=\hbox{$\vcenter{\hrule width 4pt}$}
\setbox2=\hbox{$\textstyle\int$}
\setbox\ncinttbox=\hbox{\rlap{\hbox
    to \wd2{\kern-.14em\box2\relax\hfil}}\box0\kern.1em}

\newcommand{\anticomm}[2]{\left\{ #1,#2 \right\}}
\newcommand{\cmnt}[1]{\ignorespaces}
\newcommand\norm[1]{||#1||}

\allowdisplaybreaks
\begin{document}
\title[Asymmetric noncommutative torus has vanishing Einstein tensor. manifold.]{Asymmetric noncommutative torus has vanishing Einstein tensor.} 
\author*[1]{\fnm{Deeponjit} \sur{Bose}}\email{deeponjit.bose@doctoral.uj.edu.pl}
\author[1]{\fnm{Andrzej} \sur{Sitarz}}\email{andrzej.sitarz@uj.edu.pl}
\equalcont{These authors contributed equally to this work.}
\affil*[1]{\orgdiv{Institute of Theoretical Physics}, \orgname{Jagiellonian University},  \orgaddress{\street{ul. prof. Stanis\l{}awa \L{}ojasiewicza 11}, \city{Krak\'ow}, \postcode{30-348}, \country{Poland}}}
\abstract{We explicitly compute the spectral metric, torsion and Einstein 
tensors for a nontrivial spectral triple on a noncommutative torus, with the
Dirac operator related to the fully equivariant Dirac by a partial conformal
rescaling (as introduced in \cite{DS15}). The results show that the spectral triple
has vanishing torsion and the Einstein tensor also identically vanishes.
}
\keywords{Noncommutative Geometry; spectral triple, spectral geometry, Einstein tensor.}
\pacs[MSC Classification]{58B34; 46L87}
\maketitle
\thispagestyle{empty}
\section{Introduction}
Spectral functionals, as introduced in
\cite{DSZ23} are generalization of fundamental geometric tensors in the realm of noncommutative geometry and allow to compute the metric and Einstein tensors using spectral properties of the Dirac (or Dirac type 
operators) alone. Additionally, the 
spectral methods allow to determine 
a class of torsion-free Dirac operators
without using the notion of linear connection. Though the construction has been tested on both classical and quantum geometries (including spin manifolds, Riemannian geometries with Hodge-Dirac operator, simple almost commutative geometries, conformally rescaled Dirac for the noncommutative tori and quantum $SU(2)$ group) there are still many open problems. In \cite{DSZ23} we conjectured that suitably regular two-dimensional spectral triples will have the Einstein
functional identically vanishing. 
In \cite{DS15} we computed the scalar curvature for the so-called asymmetric noncommutative torus showing that the 
integrated scalar curvature is independent of the choice of the
the Dirac operator, thus indicating this
geometry as a nice testing ground to check the validity of the conjecture.
In this paper we do prove that indeed
the Einstein functional vanishes for this geometry, supporting our conjecture.

Let us recall, that  if  $\T^2$ is the classical torus with coordinates $0 \!\leq\! x, y \!\leq\! 2\pi $, then we can equip it with the metric,
\begin{equation}
\label{metric}
dx^2 + k^{-2}(x,y) dy^2,
\end{equation}
where $k$ is a strictly positive function. This is the mildest deformation
of the flat metric. The scalar curvature of the torus with the metric 
\eqref{metric} reads 
\begin{equation}
\label{classcurv}
R = 2k^{-1} \partial_x^2 (k) - 4 k^{-2} (\partial_x (k) )^2 . 
\end{equation}
In the commutative case, such a metric is, of course, conformally equivalent to a flat metric on the torus.

Using the spectral triples approach of Connes, where the natural object is the Dirac operator rather than the metric itself, we proposed a family of Dirac operators that provide a family of  nontrivial geometries for the algebra of the noncommutative torus. In section 2, we briefly remind the construction of that spectral triple, the classical pseudodifferential calculus over the noncommutative torus, and compute the symbols of the Dirac operator and its inverse. In section 3 we recall the construction of spectral functionals and compute them for the spectral triple of the asymmetric torus. In section 4 we apply it and explicitly compute the relevant spectral functionals.

\section{The spectral triple of the asymmetric torus.}

We start with the commutative case of Dirac operator on $L^2(\T^2, k^{-1} dx\,dy)\otimes \C^2$ 
for the metric \eqref{metric}:
$$
\tilde{D}_k = -i \sigma^1 \left(\partial_x -\tfrac{1}{2} k^{-1} \partial_x (k)  \right) -  i \sigma^2 \; k\, \partial_y ,
$$
where 
\begin{equation}
\sigma^1= \left(
\begin{array}{cc}
0      &\quad 1   \\ 
 1     &\quad 0   \\ 
\end{array}
\right),\quad
\sigma^2= \left(
\begin{array}{cc}
0      &\quad -i  \\ 
 i     &\quad 0   \\ 
\end{array}
\right).
\end{equation}
It is more convenient to work with
the unitarily equivalent Dirac operator
$D_k$ acting on the dense domain of the Hilbert space $H=L^2(T,  dx\,dy)\otimes \C^2$.
\begin{equation}
\label{dirac}
D_k =  - i \sigma^1 \partial_x  -  i \sigma^2 \left(  k\, \partial_y + \tfrac{1}{2}  \partial_y(k)  \right).
\end{equation} 
The operator $D_k$ is  selfadjoint on the dense domain $H$. Next, we pass to the noncommutative torus $\T_\theta^2$, for which we refer to \cite{Co94} for the details and following \cite{DS15} we introduce the family of spectral triples.
\begin{dfn}
Let $\mathfrak{t}$ be the usual  trace on $\T_\theta^2$ and $\H= L^2(\T_\theta^2, \mathfrak{t})\otimes \C^2 $. Let $\delta_1, \delta_2$ be the usual self-adjoint derivations on $\T_\theta^2$ and  $J$ be the standard real structure
(Tomita-takesaki conjugation) on the coordinate algebra $A(\T_\theta^2)$.
For a positive element $k \in J A(\T_\theta^2) J^{-1} \subseteq  A'$, where $A'$ is the commutant of $A$, we define
\begin{equation}
\label{diracnc}
D_k =   \sigma^1 \delta_1 +   \sigma^2 \left(  k\, \delta_2 + \tfrac{1}{2}  \delta_2(k)  \right), 
\end{equation}
obtaining a family of spectral triples $(\T_\theta^2, \H,  D_k)$. 
\end{dfn}

The form of the operator $D_k$ resembles the classical Dirac operator  $D$  \eqref{dirac}, but since we take $k$ in the algebra 
$J A(\T_\theta^2) J^{-1}$ which commutes with $A(\T_\theta^2)$, 
$D_k$ has always bounded commutators with $a\in A(\T_\theta^2)$ and it gives a spectral triple in the usual sense, representing the same K-homology class
for each $k$. Moreover, $D_k$ is a differential operator 
in the sense of \cite{CoTr11} and we can use the tools of the pseudodifferential
calculus as  used in \cite{CoTr11,CoMo14}. Note that the bimodule of differential forms generated by $a [D_k,b]$ for all $a,b \in A(T_\theta^2)$ is a free left module generated by $\sigma_1$ and $k \sigma_2$. 
\section{Symbols of the Dirac}

For the details of the pseudodifferential calculus we refer the reader to \cite{CoTr11,DS15, DSZ23}.  To compute the spectral functionals,  first we have to compute the symbols of the Dirac operator and its powers in the algebra of classical pseudodifferential symbols over the noncommutative torus.
The square of $D_k$ reads,
$$
\begin{aligned}
D_k^2 &= \left( (\delta_1)^2 + k^2 (\delta_2)^2 \right), \\
          & + \left( \tfrac{3}{2} k \delta_2(k) + \tfrac{1}{2} \delta_2(k) k + i \sigma^3 \delta_1(k) \right) \delta_2, \\
          & + \left( \tfrac{1}{4}  (\delta_2(k))^2 + \tfrac{1}{2} i \sigma^3 \delta_{12}(k) + \tfrac{1}{2} k \delta_{22}(k) \right) .
\end{aligned}
$$ 
and its symbol in the pseudodifferential calculus over $A(\T_\theta^2)$ is,
$$\rho(D_k^2) = \mathfrak{a}_2 + \mathfrak{a}_1 + \mathfrak{a}_0,$$
where,
$$
\begin{aligned}
\mathfrak{a}_2 =& \left( \xi_1^2 + k^2 \xi_2^2 \right), \\
\mathfrak{a}_1 =& \left( \tfrac{3}{2} k \delta_2(k) + \tfrac{1}{2} \delta_2(k) k + i \sigma_3 \delta_1(k) \right) \xi_2, \\
\mathfrak{a}_0 =& \left( \tfrac{1}{4}  (\delta_2(k))^2 + \tfrac{1}{2} i \sigma_3 \delta_{12}(k) + \tfrac{1}{2} k \delta_{22}(k) \right).
\end{aligned}
$$
Next, we compute the first three leading symbols of $D_k^{-1}$,
\begin{align}
\rho(D_k^{-1}) = \mathfrak{b}_{-1} + \mathfrak{b}_{-2} + \mathfrak{b}_{-3} + \dots ,
\end{align}
where $b_{-n}$ is a  symbol of $D_k^{-1}$ homogeneous of degree $(-n) $ in $\xi$. We have
$$ \mathfrak{b}_{-1} = (\sigma^1 \xi_1 + k \sigma^2 \xi_2) (\xi_1^2 + k^2 \xi_2^2)^{-1}. $$
The next homogeneous symbols we obtain from:
$$ \mathfrak{b}_{1} \mathfrak{b}_{-2}   +  \mathfrak{b}_{0} \mathfrak{b}_{-1} 
+ \partial_{\xi^j} (\mathfrak{b}_{1}) \delta_j (\mathfrak{b}_{-1}) = 0 , $$
which gives:
$$ 
\mathfrak{b}_{-2} = - \mathfrak{b}_{-1}  \bigl( \mathfrak{b}_{0} \mathfrak{b}_{-1} + \partial_{\xi^j} (\mathfrak{b}_{1}) \delta_{j} (\mathfrak{b}_{-1}) \bigr),
$$
and that can be computed explicitly:
\begin{equation}
\begin{aligned}
\mathfrak{b}_{-2} = & \sigma^1  \biggl( k b_0 \delta_1(k) b_0 \xi_2^2 - 
\frac{3}{2}  k b_0 \delta_2(k) 
         b_0 \xi_1 \xi_2 + k b_0^2 \delta_1(k) b_0 \xi_1^2 \xi_2^2 \\ 
         &\quad- k^2 b_0^2
          \delta_1(k) k b_0 \xi_2^4 + 2 k^2 b_0^2 \delta_2(k) k b_0 \xi_1 
         \xi_2^3 - k^3 b_0^2 \delta_1(k) b_0 \xi_2^4 \\
         &\quad + 2 k^3 b_0^2 \delta_2(k) 
          b_0 \xi_1 \xi_2^3 
         - \frac{1}{2}  b_0 \delta_2(k) k b_0 \xi_1 \xi_2 + b_0^2 
         \delta_1(k) k b_0 \xi_1^2 \xi_2^2 \biggr) \\
&       + \sigma^2  \biggl(  - \frac{1}{2} k b_0 \delta_2(k) k b_0 \xi_2^2 + 2 k b_0^2 
         \delta_1(k) k b_0 \xi_1 \xi_2^3 - k b_0^2 \delta_2(k) k b_0 \xi_1^2 
         \xi_2^2 \\
         & \quad - k^2 b_0 \delta_2(k) b_0 \xi_2^2 + 2 k^2 b_0^2 \delta_1(k) 
         b_0 \xi_1 \xi_2^3 - k^2 b_0^2 \delta_2(k) b_0 \xi_1^2 \xi_2^2 + k^3 
         b_0^2 \delta_2(k) k b_0 \xi_2^4 \\
         &\quad+ k^4 b_0^2 \delta_2(k) b_0 \xi_2^4
          - b_0 \delta_1(k) b_0 \xi_1 \xi_2 + \frac{1}{2} b_0 \delta_2(k) b_0 \xi_1^2
          \biggr),
\end{aligned}
\label{b-2}
\end{equation}
where by $b_0$ we denote $b_0 = (\xi_1^2 + k^2 \xi_2^2)^{-1}$.

In the next step  we compute the three leading symbols of $D_k^{-2}$, 
\begin{align}
  \rho(D_k^{-2}) = \mathfrak{c}_{-2} + \mathfrak{c}_{-3} + \mathfrak{c}_{-4} + \ldots
\end{align}
where again, $\mathfrak{c}_{-n}$ is a homogeneous symbol of $D_k^2$ of degree $(-n)$ in $\xi$. 
We obtain,

\begin{subequations}
    \begin{align}
        \mathfrak{c}_{-2} &= (\mathfrak{a}_2)^{-1} = (\xi_1^2 + k^2 \xi_2^2)^{-1} = b_0,
        \\
        \nonumber
        \mathfrak{c}_{-3} &= -\mathfrak{c}_{-2} \left(\mathfrak{a}_1 \mathfrak{c}_{-2} + \sum_{i=1, 2}\left(\del_{\xi^i}(\mathfrak{a}_2)\, (\delta_i  (\mathfrak{c}_{-2})\right)\right)
        \\ 
        &= 
        \begin{aligned}[t] 
            &
            2 \left(\xi _2\right)^3 (\mathfrak{b}_{-2})^2 k^3\delta _2(k)\mathfrak{b}_{-2}+2 \left(\xi _2\right){}^3 (\mathfrak{b}_{-2})^2 k^2\delta _2(k)k \mathfrak{b}_{-2} - i \sigma^3 \xi_2 \mathfrak{b}_{-2}\,\delta _1(k) \,\mathfrak{b}_{-2}
            \\
            & \hspace{0.2cm}
            + 2 \, \xi _1 \left(\xi _2\right)^2 (\mathfrak{b}_{-2})^2 \left(k\delta _1(k)
            + \delta _1(k) k \right)\mathfrak{b}_{-2}-\frac{3}{2} \xi _2 \mathfrak{b}_{-2} k\delta _2(k)\mathfrak{b}_{-2}
            \\
            & \hspace{0.2cm}
            -\frac{1}{2} \xi_2 \mathfrak{b}_{-2} \delta _2(k)k \mathfrak{b}_{-2},
\end{aligned}
        \\
        \mathfrak{c}_{-4} &= \begin{aligned}[t]
            & - \mathfrak{c}_2\left(\mathfrak{a}_1\mathfrak{c}_{-3} + \mathfrak{a}_0\mathfrak{c_{-2}} + \sum_{j=1, 2} \left(\del_{\xi^j}(\mathfrak{a}_1)\delta_j(\mathfrak{c}_{-2}) + \del_{\xi^j}(\mathfrak{a}_2)\delta_j(\mathfrak{c}_{-3})\right)\right.
            \\
            &\hspace{2cm}
            \left. +  \frac{1}{2}\sum_{l, m = 1 ,2}\left(\del_l \del_m (\mathfrak{a}_2) \, \delta_l \delta_m (\mathfrak{b}_{-2})\right)\right).
        \end{aligned}
    \end{align}\label{symbolsDinv}
\end{subequations}
\section{Spectral functionals}
In \cite{DSZ23, DSZ24,DSZ24a} a family of functionals was {introduced} and
studied for finitely summable spectral triples which have a Wodzicki-residue noncommutative trace \cite{Wo87} (denoted as $\operatorname{Wres}$) over the Connes-Moscovici pseudodifferential calculus.
In the classical case of spectral triples over smooth functions of a Riemannian 
spin manifold and the Dirac operator being the usual spin Dirac differential operator these functional were demonstrated to reproduce the evaluation of geometric tensors
on differential forms. The construction was shown to give the same tensors also for 
the Hodge-Dirac operator and was succesfully applied to some genuine noncommutative examples. 

We shall now compute explicitly these functional for the spectral triple defined
in the section 2.

\subsection{The metric functional}
Let us recall, that for any two one-forms in $\Omega^1_D(A)$ the metric
functional for a $n$-summable spectral triple over $A$ is,
\begin{equation}
	\cg_{D}(u, v) = \operatorname{Wres}\left(u\, v \, |D|^{-n}\right) .
\end{equation} 
\begin{pro}
Let $u_i,v_j \in A(\T^2_\theta)$, $i,j=1,2$ and $u = u_1 \sigma_1 + k  u_2 \sigma_2$,
$v = v_1 \sigma_1 + k  v_2 \sigma_2$ be the one-forms in $\Omega^1(A(\T^2_\theta))$. 
Then the metric functional for the spectral triple of the asymmetric noncommutative
torus is:	
\begin{equation}
	\cg_{D_k}(u, v) = \tau\left(\frac{1}{k}u_1 v_1 + k u_2 v_2 \right)
\end{equation} 
\end{pro}
\noindent{\sl Proof:~~~} 
We compute using results from section 3:
\begin{align}
	\cg_{D_k}(u, v) &= \operatorname{Wres}\left(u\, v \, D_k^{-2}\right) = \int_{\norm{\xi}=1} \mathfrak{t} \left(\operatorname{Tr}\left(\rho_{-2}\left(u\, v \, D_{k}^{-2}\right)\right)\right)d\xi \\
	& =  \int_{\norm{\xi}=1} \mathfrak{t}
	\biggl( \left( u_1 v_1 + k^2 u_2 v_2 \right) (\xi_1^2 + k^2 \xi_2^2)^{-1}  \biggr) \\
	& = 2\pi \tau\left(\frac{1}{k}u_1 v_1 + k u_2 v_2 \right).
\end{align}
where we have use that trace over the product $\sigma_1 \sigma_2$ 
vanishes and $k$ commutes with $u_i,v_j$.
 $\qquad\Box$
\subsection{The torsion functional}
The torsion functional is a functional over three one-forms $u,v,w$, which
is defined, for an $n$-dimensional spectral geometry, as
\begin{equation}
	{\mathcal T}_{D}(u, v, w) = 
	\operatorname{Wres}\left(u\, v \, w \, D |D|^{-n}\right) ,
\end{equation} 
and, in the classical case it detects whether the Dirac operator $D$ has torsion,
that is, whether it is a lift to spinors of a linear connection with a nonvanishing 
torsion. 
In our case, to facilitate the computation of  the torsion functional, observe that the product of three one-forms will necessarily be a linear combination of the form $a \sigma^1 + b\sigma^2$, where  $a,b \in A(\T^2_\theta)$ and they commute 
with $b_0$ and $k$.
 
\begin{pro}
The torsion functional for the spectral triple of the asymmetric noncommutative
torus vanishes identically.
\end{pro} 
\noindent{\sl Proof:~~~} 
Using the argument for the product of $uvw$ it is sufficient that we  compute
for arbitrary $a,b \in A(\T^2_\theta)$ :
\begin{equation}
 \begin{aligned}
 \mathcal{W} &\biggl( (a \sigma^1 + b\sigma^2) D_k |D_k|^{-2} \biggr)= 
 \int_{S^1} \tau \biggl((a \sigma^1 + b\sigma^2 ) \mathfrak{b}_{-2} \biggr) \\
  &= \int_{S^1}  \tau \biggl(a \biggl(
  (  - 2\delta_1(k) k^3 b_0^3 \xi_2^4 + 
         \delta_1(k) k b_0^2 \xi_2^2 + 2 \delta_1(k) k b_0^3 \xi_1^2 \xi_2^2
         \biggr) \biggr) \\
     &+    \int_{S^1}  \tau \biggl(b \biggl( -\frac{3}{2}  \delta_2(k) k^3 b_0^2 \xi_2^2 
        + \frac{1}{2} \delta_2(k) k b_0^2 \xi_1^2  
         - 2\delta_2(k) k^3 b_0^3 \xi_1^2 \xi_2^2 
         + 2\delta_2(k) k^4 b_0^3 \xi_2^4 \biggr) \biggr) 
 \end{aligned} 	
\end{equation}
Computing the integrals explicitly we have:
\begin{equation}
\begin{aligned}
	\mathcal{W} \biggl( (a \sigma^1 + b\sigma^2) D_k |D_k|^{-2} \biggr) 
	& = \tau \biggl( a \delta_1(k) \left( -\frac{3\pi}{2k^2}
    + \frac{\pi}{k^2} + \frac{\pi}{2k^2} \right) \biggr) 
    \\
   & \hspace{-0.5cm} + \tau \biggl( b \delta_2(k) \left(-\frac{3}{2} \pi
  + \frac{1}{2} \pi - 2 \frac{\pi}{4}
     + 2\frac{3\pi}{4}\right) \biggr) = 0.
\end{aligned}
\end{equation}
 $\Box$

\noindent Please observe that although the classical geometry of the torus is
two-dimensional, hence  there is no fully antisymmetric torsion and therefore
all Dirac operators have no torsion. However, this argument does not work
in the noncommutative situation and only explicit computations allow us to
demonstrate this feature of the spectral triple and the Dirac operator $D_k$.

Let us finish this part by observing that the spectral triple with the operator $D_k$ is, in fact,
spactrally closed in the sense of  \cite{DSZ23}, that is for $T$, which is a product of an arbitrary
numbers of one-forms, the residue:
$$ \mathcal{W} \biggl( T D_k |D_k|^{-2} \biggr)= 0. $$
The proof is similar to the above for torsion.
\section{The spectral Einstein functional}
The spectral Einstein function a pair of 1-forms $u, v$ is given as
\begin{align}
	\mathcal{G}_{D}(u,v) = \operatorname{Wres}\left(u\, \anticomm{D}{v} \, D \,D^{-n}\right).
\end{align}
for an $n$-dimensional spectral geometry.

\begin{align}
    \mathcal{G}_{D_k} = \operatorname{Wres}\left(u\, \anticomm{D_k}{v} \, D_k \,D_k^{-2}\right) = \int_{\norm{\xi}=1} \mathfrak{t} \left(\operatorname{Tr}\left(\rho_{-2}\left(u\, \anticomm{D_k}{v} \, D_k^{-1}\right)\right)\right)d\xi
\end{align}

To compute the symbol $\rho_{-2}\left(u\, \anticomm{D_k}{v} \, D_k^{-1}\right)$ we first compute 
the product $u \{ D_k, v\}$ for any differential forms $u,v$ where 
$u= u_1 \sigma_1 + u_2 k\sigma_2$ and $v= v_1 \sigma_1 + v_2 k\sigma_2$:
\begin{align}
    u\,\anticomm{D_k}{v} = 
    \begin{aligned}[t]
        &
        2\sigma ^2k^3u_2v_2\,\delta _2 + k^3 \sigma ^2 u_2\delta _2\left(v_2\right) + \frac{3}{2} \sigma^2 k^2 \delta _2(k) u_2 v_2
        \\
        & \hspace{0.2cm}
        +2 \sigma ^1k^2 u_1v_2 \, \delta _2 + \sigma^1 k^2 u_1 \delta _2\left(v_2\right) -\sigma^1 k^2  u_2 \delta _1\left(v_2\right) 
        \\
        & \hspace{0.2cm}
        + \sigma ^2 \delta _1(k) u_1v_2 + \sigma^2 k u_1\delta _1\left(v_2\right) + 2 \sigma ^2 k u_2 v_1\, \delta _1 + k\sigma ^2 u_2\delta _1\left(v_1\right)
        \\
        & \hspace{0.2cm}
        +\frac{1}{2} \sigma ^2 k\delta _2(k)k u_2 v_2 +\frac{3}{2} \sigma ^1 k\delta _2(k) u_1 v_2 +\frac{1}{2} \sigma ^1 \delta _2(k)k u_1v_2 
        \\
        & \hspace{0.2cm}
          + \sigma ^1 k^2 u_2 \delta _2\left(v_1\right) - \sigma ^1 k\delta _1(k) u_2v_2 - \sigma ^2 k u_1\delta _2\left(v_1\right)
        \\
        & \hspace{0.2cm}
         + 2 \sigma ^1 u_1v_1\, \delta _1 + \sigma ^1 u_1\delta _1\left(v_1\right)
    \end{aligned}
\end{align}
This is a first order differential operator on $A(T_\theta^2)$, so its homogeneous symbols of order 1 and 0 are,
\begin{align}
    \rho_{1}\left(u\anticomm{D_k}{v}\right)
    &= 2 \, \left(\xi_2 \, u_1v_2k^2 \sigma ^1 + \xi _1 u_2v_1k\sigma ^2 + \xi _1 u_1v_1\sigma ^1 + \xi _2 u_2v_2 k^3 \sigma ^2\right)
    \\
    \rho_{0}\left(u\anticomm{D_k}{v}\right)
    &=
    \begin{aligned}[t]
        &
        u_1\delta _2\left(v_2\right)\sigma ^1k^2+u_2\delta _2\left(v_1\right)\sigma ^1k^2+u_1v_2\delta _1(k)\sigma ^2+u_1\delta _1\left(v_2\right)k\sigma ^2+u_2\delta _1\left(v_1\right)k\sigma ^2
        \\
        & \hspace{0.2cm}
        +\frac{3}{2} u_1v_2k\delta _2(k)\sigma ^1+\frac{1}{2} u_1v_2\delta _2(k)k\sigma ^1-u_1\delta _2\left(v_1\right)k\sigma ^2+u_1\delta _1\left(v_1\right)\sigma ^1
        \\
        & \hspace{0.2cm}
        +\sigma ^2 u_2\delta _2\left(v_2\right)k^3
        +\frac{3}{2} u_2v_2k^2\delta _2(k)\sigma ^2-u_2\delta _1\left(v_2\right)\sigma ^1k^2+\frac{1}{2} u_2v_2k\delta _2(k)k\sigma ^2
        \\
        & \hspace{0.2cm}
        -u_2v_2k\delta _1(k)\sigma ^1
    \end{aligned}
\end{align}
Finally, we can compute the symbol of order $-2$ of  $u \anticomm{D_k}{v} D_k^{-1}$:
\begin{equation}
\begin{aligned}
    \rho_{-2}  \left(u \anticomm{D_k}{v} D_k^{-1}\right) = &
        \rho_{0}\left(u\anticomm{D_k}{v}\right)\, \mathfrak{b}_{-2} 
        + \rho_{1}\left(u\anticomm{D_k}{v}\right)\, \mathfrak{b}_{-3}
        \\
        &\hspace{0.8cm}
        + \, \del_{\xi^1} \left(\rho_{1}\left(u\anticomm{D_k}{v}\right)\right) \delta_1 \mathfrak{b}_{-2}  \\
        &\hspace{0.8cm}
        + \, \del_{\xi^2} \left(\rho_{1}\left(u\anticomm{D_k}{v}\right)\right) \delta_2 \mathfrak{b}_{-3}    \\
        &\hspace{0.8cm}
         +\del_{\xi^1}\del_{\xi^2}  \left(\rho_{1}\left(u\anticomm{D_k}{v}\right)\right) \delta_1\delta_2 \left( \mathfrak{b}_{-1}\right).
\label{symbein}
\end{aligned}
\end{equation}
The explicit form of the symbol \eqref{symbein}, using the computed symbols of $D_k^{-1}$ \eqref{b-2}, is provided explicitly in the Appendix.

\subsection{Evaluating the Wodzicki residue.}

The major task in the computation of the Wodzicki residue of  \eqref{symbein}
is the calculation of the integral of this expression over $S^1$ in the $\{\xi_1,\xi_2\}$ plane. First, it is easily verified that the sum of the integrals of the coefficients $u_i\delta_j(v_k)$ vanish identically. This is in total agreement with observation
in \cite{DSZ23} that if the Dirac operator has no torsion then the Einstein
functional is tensorial and thus does not depends on derivations of the 
coefficients of one forms and is intrinsically bilinear, that is,
\begin{align}
	\mathcal{G}_{D}(au,vb) = \operatorname{Wres}\left(au\, \anticomm{D}{v}b \, D \,D^{-n}\right).
\end{align}
Therefore we can focus on the part of the symbol of $\rho_{-2}(u\anticomm{D_k}{v})$  excluding all such terms. The resulting expression consists of 320 terms in the form of products $u_i v_j$ multiplied by the a function of $k$ and its derivates.  Next, we check that the sum the terms, which contain the second derivatives of $k$, $\delta_i\delta_j(k)$ also vanish.  All remaining terms are of the form,  for some $m,n,a,b,\alpha,\beta$:
$$
I(m,n,a,b,\alpha,\beta) = \int_0^{2\pi} d\phi \, \tau \biggl( b_0^m k^a \delta_i(k) b_0^n k^b \delta_j(k) \xi_1^{2\alpha} \xi_2^{2\beta} \biggr), 
$$
where $\xi_1 = \cos\phi$ and $\xi_2=\sin\phi$. Using a simple change of variables:
$$ z = \tan(\phi), $$
it can be rewritten as ,
\begin{align}
    I(m,n,a,b,\alpha,\beta) = 4 \int_0^\infty dz \, 
(1+z^2)^{m+n-\alpha-\beta-1} 
\tau \biggl( k^a \biggl( \frac{1}{1+k^2z^2} \biggr)^m
\delta_i(k) k^b  \biggl( \frac{1}{1+k^2z^2} \biggr)^n 
\delta_j(k) z^{2\beta} \biggr). 
\end{align}
At this point we can use (a slightly modified) Lesch rearrangement lemma
\cite{Le17, HeMDvN24}. Introducing the operator $\Delta$ acting on the elements of $A(\T^2_\theta)$, {with $k$ being a positive operator}:
$$\Delta(A) = k^{-1} A k,$$
we can express $I$ as:
\begin{equation}
\begin{aligned}
	&I(m,n,a,b,\alpha,\beta) = \\
	&= 4  \tau \biggl( \biggl[ \int_0^\infty dz \, 
	(1+z^2)^{m+n-\alpha-\beta-1} z^{2\beta} k^a \biggl( \frac{1}{1+k^2z^2} \biggr)^m k^b \Delta^b \biggl( \frac{1}{1+k^2z^2 \Delta^2} \biggr)^n \biggr] \biggl( \delta_i(k) \biggr) \cdot  \delta_j(k) \biggr),
\end{aligned}
\end{equation}
where the operator in square brackets acts on $\delta_i(k)$. 
Now, since each expression  is a part of a symbol homogeneous in $\xi$ of order -2,  the sum $m + n-\alpha -\beta -1$ has to vanish, so the the first term $(1+z^2)^{m+n-\alpha-\beta-1}$ is not present, which allows us to change  variable $u=kz$ and express $I$ as,
\begin{align}
    I(m,n,a,b,\alpha,\beta) =  \tau \biggl( \bigl[  {F}(\Delta; m,n,\beta, a,b ) (\delta_i(k))  \bigr] \delta_j(k) \biggr),
\end{align}    
 where   
\begin{align}
{F}(s; m,n,\beta,a,b) =   4  k^{a+b-2\beta-1}  \int_0^\infty du \, u^{2\beta} 
\biggl( \frac{1}{1+u^2} \biggr)^m s^b 
 \biggl( \frac{1}{1+u^2 s^2} \biggr)^n.
\end{align}
We explicitly compute all 290 integrals and show  all results  in tables (\ref{tab1}-\ref{tab6}) in the Appendix.

\section{Conclusions}
The {noncommutative} tori, and in particular the two-dimensional noncommutative torus,  are a nice testing ground for noncommutative geometries. The conformally rescaled geometry of the noncommutative torus has been the subject of many papers \cite{CoTr11, FaKh12, FaKh13, FaKh13a, Ro13, CoMo14} showing that Gauss-Bonnet
theorem holds. The approach, though not based on the usual spectral triple, to study the Ricci curvature was proposed in \cite{FGK19} and the Einstein tensor was showed to vanish in \cite{DSZ23}. The asymmetric noncommutative torus, which was proposed and studied in \cite{DS15} is yet another incarnation of  an interesting family of geometries for this noncommutative manifold. It is highly interesting, because the only way to prove that the Gauss-Bonnet theorem holds and the Einstein tensor vanishes was by an explicit computation.

The results we obtained here confirm that the spectral triple with the partial conformal rescaling of the Dirac operator is a valid, two-dimensional geometry. First of all, the Dirac has vanishing torsion (in fact, the spectral triple is spectrally closed), which is not guranteed in noncommutative geometry. Then, the vanishing of the Einstein functional is consistent with the Gauss-Bonnet theorem, which holds for this geometry and shows
that the construction of Einstein functional has a deep geometrical sense not only in the classical sense. Although this is just a single probe into the
spectral geometry of a broad family of Dirac operators over the noncommutative tori, it confirms that at least some two-dimensional noncommutative geometries behave similarly to the classical ones.

\vspace{1cm}

{\bf Acknowledgements} D.B. gratefully acknowledges the hospitality of SISSA and Ludwik D\k{a}browski as well as support by research support mudule RSM/76/BD
from the Jagiellonian University. 

This is a part of the international project  Graph Algebras partially supported by EU grant HORIZON-MSCA SE-2021 Project 101086394 and co-financed by the Polish Ministry of Education 
and Science within the program PMW under contract 5447/HE/2023/2 


\appendix
\pagebreak[4]
\section{Appendix: the symbol.}

First, we write explicitly the form of the symbol \eqref{symbein} that after computing the Wodzicki residue gives the Einstein functional. {We write down only the even powers of $\xi_1, \xi_2$ since the odd powers trivially vanish under the integral over $S^1$ in $\xi_1 - \xi_2$ plane:}
\begin{align*}
\tau\Big(\rho^{\text{even}}_{-2}  \bigl(u \anticomm{D_k}{v} & D_k^{-1} \bigr)\Big) \\
=	{\tau\Bigg(}u_1 v_1 \Big(
		& \left(12 \left(\xi _2\right){}^4 \left(\xi _1\right){}^2 b_0^4k^3-4 \left(\xi _2\right){}^4 b_0^3k^3-4 \left(\xi _2\right){}^2 \left(\xi _1\right){}^4 b_0^4k+2 \left(\xi _2\right){}^2 b_0^2k\right)\delta _1\left(\delta _1(k)\right) 
		\\
		& + \left(-12 \left(\xi _2\right){}^4 \left(\xi _1\right){}^2 b_0^4k^5+4 \left(\xi _2\right){}^2 \left(\xi _1\right){}^4 b_0^4k^3+7 \left(\xi _2\right){}^2 \left(\xi _1\right){}^2 b_0^3k^3-\left(\xi _1\right){}^4 b_0^3k\right)\delta _2\left(\delta _2(k)\right)
		\\
		& +  \Big(6 b_0^3k^4\delta _1(k)b_0 \left(\xi _2\right){}^6-2 \left(\xi _1\right){}^2 b_0^3k^4\delta _1(k)b_0^2 \left(\xi _2\right){}^6-18 \left(\xi _1\right){}^2 b_0^4k^4\delta _1(k)b_0 \left(\xi _2\right){}^6
		\\
		& + 2 b_0^2k\delta _1(k)b_0^2k^3 \left(\xi _2\right){}^6+4 b_0^2k^2\delta _1(k)b_0^2k^2 \left(\xi _2\right){}^6+2 b_0^2k^3\delta _1(k)b_0^2k \left(\xi _2\right){}^6
		\\
		& -4 \left(\xi _1\right){}^2 b_0^3k\delta _1(k)b_0^2k^3 \left(\xi _2\right){}^6+2 b_0^3k^2\delta _1(k)b_0k^2 \left(\xi _2\right){}^6-10 \left(\xi _1\right){}^2 b_0^3k^2\delta _1(k)b_0^2k^2 \left(\xi _2\right){}^6
		\\
		& + 8 b_0^3k^3\delta _1(k)b_0k \left(\xi _2\right){}^6-8 \left(\xi _1\right){}^2 b_0^3k^3\delta _1(k)b_0^2k \left(\xi _2\right){}^6-6 \left(\xi _1\right){}^2 b_0^4k^2\delta _1(k)b_0k^2 \left(\xi _2\right){}^6
		\\
		& -24 \left(\xi _1\right){}^2 b_0^4k^3\delta _1(k)b_0k \left(\xi _2\right){}^6-2 b_0^3k^2\delta _1(k){}^2 \left(\xi _2\right){}^4+6 \left(\xi _1\right){}^2 b_0^4k^2\delta _1(k){}^2 \left(\xi _2\right){}^4
		\\
		&-2 b_0\delta _1(k)b_0^2k^2 \left(\xi _2\right){}^4-8 b_0^2k^2\delta _1(k)b_0\left(\xi _2\right){}^4+12 \left(\xi _1\right){}^2 b_0^3k^2\delta _1(k)b_0 \left(\xi _2\right){}^4
		\\
		& +2 \left(\xi _1\right){}^4 b_0^3k^2\delta _1(k)b_0^2\left(\xi _2\right){}^4-2 \left(\xi _1\right){}^2 b_0^3\delta _1(k)b_0k^2 \left(\xi _2\right){}^4+2 \left(\xi _1\right){}^4 b_0^3\delta _1(k)b_0^2k^2 \left(\xi _2\right){}^4
		\\
		& +6 \left(\xi _1\right){}^4 b_0^4k^2\delta _1(k)b_0 \left(\xi _2\right){}^4+2 \left(\xi _1\right){}^4 b_0^4\delta _1(k)b_0k^2 \left(\xi _2\right){}^4-2 b_0k\delta _1(k)b_0^2k \left(\xi _2\right){}^4
		\\
		& -6 b_0^2k\delta _1(k)b_0k\left(\xi _2\right){}^4+4 \left(\xi _1\right){}^2 b_0^3k\delta _1(k)b_0k\left(\xi _2\right){}^4+4 \left(\xi _1\right){}^4 b_0^3k\delta _1(k)b_0^2k\left(\xi _2\right){}^4
		\\
		& +8 \left(\xi _1\right){}^4 b_0^4k\delta _1(k)b_0k\left(\xi _2\right){}^4+2 \left(\xi _1\right){}^2 b_0^3\delta _1(k) \left(\xi _2\right){}^2-2 \left(\xi _1\right){}^4 b_0^4\delta _1(k) \left(\xi _2\right){}^2
		\\
		& +2 b_0\delta _1(k)b_0\left(\xi _2\right){}^2-2 \left(\xi _1\right){}^4 b_0^3\delta _1(k)b_0 \left(\xi _2\right){}^2 \Big)\delta _1(k) + \Big( 18 \left(\xi _1\right){}^2 b_0^4k^6\delta _2(k)b_0\left(\xi _2\right){}^6
		\\
		& +4 \left(\xi _1\right){}^2 b_0^3k^2\delta _2(k)b_0^2k^4 \left(\xi _2\right){}^6+10 \left(\xi _1\right){}^2 b_0^3k^3\delta _2(k)b_0^2k^3 \left(\xi _2\right){}^6
		\\
		& + 8 \left(\xi _1\right){}^2 b_0^3k^4\delta _2(k)b_0^2k^2 \left(\xi _2\right){}^6 
		+2 \left(\xi _1\right){}^2 b_0^3k^5\delta _2(k)b_0^2k \left(\xi _2\right){}^6
		\\
		&+6 \left(\xi _1\right){}^2 b_0^4k^4\delta _2(k)b_0k^2 \left(\xi _2\right){}^6+24 \left(\xi _1\right){}^2 b_0^4k^5\delta _2(k)b_0k \left(\xi _2\right){}^6-6 \left(\xi _1\right){}^2 b_0^4k^4\delta _2(k)\left(\xi _2\right)^4
		\\
		&-\left(\xi _1\right){}^2 b_0^2\delta _2(k)b_0^2k^4 \left(\xi _2\right){}^4-20 \left(\xi _1\right){}^2 b_0^3k^4\delta _2(k)b_0 \left(\xi _2\right){}^4-6 \left(\xi _1\right){}^4 b_0^4k^4\delta _2(k)b_0\left(\xi _2\right)^4
		\\
		& -5 \left(\xi _1\right){}^2 b_0^2k\delta _2(k)b_0^2k^3 \left(\xi _2\right){}^4-6 \left(\xi _1\right){}^2 b_0^2k^2\delta _2(k)b_0^2k^2\left(\xi _2\right){}^4
		\\
		&-2 \left(\xi _1\right){}^2 b_0^2k^3\delta _2(k)b_0^2k \left(\xi _2\right){}^4
		-2 \left(\xi _1\right){}^4 b_0^3k\delta _2(k)b_0^2k^3 \left(\xi _2\right)^4
		\\
		& -10 \left(\xi _1\right){}^2 b_0^3k^2\delta _2(k)b_0k^2 \left(\xi _2\right){}^4-4 \left(\xi _1\right){}^4 b_0^3k^2\delta _2(k)b_0^2k^2 \left(\xi _2\right){}^4
		\\
		&-24 \left(\xi _1\right){}^2 b_0^3k^3\delta _2(k)b_0k \left(\xi _2\right){}^4-2 \left(\xi _1\right){}^4 b_0^3k^3\delta _2(k)b_0^2k \left(\xi _2\right)^4
		\\
		&-2 \left(\xi _1\right){}^4 b_0^4k^2\delta _2(k)b_0k^2 \left(\xi _2\right){}^4-8 \left(\xi _1\right){}^4 b_0^4k^3\delta _2(k)b_0k \left(\xi _2\right){}^4+6 \left(\xi _1\right){}^2 b_0^3k^2\delta _2(k) \left(\xi _2\right){}^2
		\\
		&+2 \left(\xi _1\right){}^4 b_0^4k^2\delta _2(k) \left(\xi _2\right){}^2+\frac{5}{2} \left(\xi _1\right){}^2 b_0^2k^2\delta _2(k)b_0 \left(\xi _2\right){}^2
		+\left(\xi _1\right){}^2 b_0^2\delta _2(k)b_0k^2 \left(\xi _2\right){}^2
		\\
		&+\left(\xi _1\right){}^4 b_0^2\delta _2(k)b_0^2k^2 \left(\xi _2\right){}^2+6 \left(\xi _1\right){}^4 b_0^3k^2\delta _2(k)b_0 \left(\xi _2\right){}^2+4 \left(\xi _1\right){}^2 b_0^2k\delta _2(k)b_0k\left(\xi _2\right){}^2
		\\
		&+\left(\xi _1\right){}^4 b_0^2k\delta _2(k)b_0^2k\left(\xi _2\right){}^2+4 \left(\xi _1\right){}^4 b_0^3k\delta _2(k)b_0k\left(\xi _2\right){}^2-\frac{1}{2} \left(\xi _1\right){}^4 b_0^2\delta _2(k)b_0 \Big)\delta _2(k) \Big) \\
	 + u_1 v_2\Big(
		& \left(8 \left(\xi _2\right){}^6 b_0^4k^7-9 \left(\xi _2\right){}^4 b_0^3k^5-24 \left(\xi _1\right){}^2 \left(\xi _2\right){}^4 b_0^4k^5+2 \left(\xi _2\right){}^2 b_0^2k^3+11 \left(\xi _1\right){}^2 \left(\xi _2\right){}^2 b_0^3k^3\right)\delta _1\left(\delta _2(k)\right)
		\\
		& +\Big(-6 b_0^4k^8\delta _1(k)b_0\delta _2(k) \left(\xi _2\right){}^8-6 b_0^4k^8\delta _2(k)b_0\delta _1(k) \left(\xi _2\right){}^8-2 b_0^3k^4\delta _1(k)b_0^2k^4\delta _2(k) \left(\xi _2\right){}^8
		\\
		&-4 b_0^3k^5\delta _1(k)b_0^2k^3\delta _2(k) \left(\xi _2\right){}^8-2 b_0^3k^5\delta _2(k)b_0^2k^3\delta _1(k) \left(\xi _2\right){}^8-2 b_0^3k^6\delta _1(k)b_0^2k^2\delta _2(k) \left(\xi _2\right){}^8		
		\\
		&-4 b_0^3k^6\delta _2(k)b_0^2k^2\delta _1(k) \left(\xi _2\right){}^8-2 b_0^3k^7\delta _2(k)b_0^2k\delta _1(k) \left(\xi _2\right){}^8-2 b_0^4k^6\delta _1(k)b_0k^2\delta _2(k) \left(\xi _2\right){}^8
		\\
		&-2 b_0^4k^6\delta _2(k)b_0k^2\delta _1(k) \left(\xi _2\right){}^8-8 b_0^4k^7\delta _1(k)b_0k\delta _2(k) \left(\xi _2\right){}^8-8 b_0^4k^7\delta _2(k)b_0k\delta _1(k) \left(\xi _2\right){}^8
		\\
		&+2 b_0^4k^6\delta _1(k)\delta _2(k) \left(\xi _2\right){}^6
		+2 b_0^4k^6\delta _2(k)\delta _1(k) \left(\xi _2\right){}^6+11 b_0^3k^6\delta _1(k)b_0\delta _2(k) \left(\xi _2\right){}^6
		\\
		&+10 b_0^3k^6\delta _2(k)b_0\delta _1(k) \left(\xi _2\right){}^6+4 \left(\xi _1\right){}^2 b_0^3k^6\delta _2(k)b_0^2\delta _1(k) \left(\xi _2\right){}^6
		+18 \left(\xi _1\right){}^2 b_0^4k^6\delta _1(k)b_0\delta _2(k) \left(\xi _2\right){}^6
		\\
		&+18 \left(\xi _1\right){}^2 b_0^4k^6\delta _2(k)b_0\delta _1(k) \left(\xi _2\right){}^6+2 b_0^2k^3\delta _1(k)b_0^2k^3\delta _2(k) \left(\xi _2\right){}^6+3 b_0^2k^3\delta _2(k)b_0^2k^3\delta _1(k) \left(\xi _2\right){}^6
		\\
		&+2 b_0^2k^4\delta _1(k)b_0^2k^2\delta _2(k) \left(\xi _2\right){}^6+7 b_0^2k^4\delta _2(k)b_0^2k^2\delta _1(k) \left(\xi _2\right){}^6+4 b_0^2k^5\delta _2(k)b_0^2k\delta _1(k) \left(\xi _2\right){}^6
		\\
		&+2 \left(\xi _1\right){}^2 b_0^3k^2\delta _1(k)b_0^2k^4\delta _2(k) \left(\xi _2\right){}^6+8 \left(\xi _1\right){}^2 b_0^3k^3\delta _1(k)b_0^2k^3\delta _2(k) \left(\xi _2\right){}^6		
		\\
		&+2 \left(\xi _1\right){}^2 b_0^3k^3\delta _2(k)b_0^2k^3\delta _1(k) \left(\xi _2\right){}^6+4 b_0^3k^4\delta _1(k)b_0k^2\delta _2(k) \left(\xi _2\right){}^6+10 \left(\xi _1\right){}^2 b_0^3k^4\delta _1(k)b_0^2k^2\delta _2(k) \left(\xi _2\right){}^6
		\\
		& + 3 b_0^3k^4\delta _2(k)b_0k^2\delta _1(k) \left(\xi _2\right){}^6+8 \left(\xi _1\right){}^2 b_0^3k^4\delta _2(k)b_0^2k^2\delta _1(k) \left(\xi _2\right){}^6+13 b_0^3k^5\delta _1(k)b_0k\delta _2(k) \left(\xi _2\right){}^6
		\\
		&+4 \left(\xi _1\right){}^2 b_0^3k^5\delta _1(k)b_0^2k\delta _2(k) \left(\xi _2\right){}^6+11 b_0^3k^5\delta _2(k)b_0k\delta _1(k) \left(\xi _2\right){}^6+10 \left(\xi _1\right){}^2 b_0^3k^5\delta _2(k)b_0^2k\delta _1(k) \left(\xi _2\right){}^6
		\\
		&+6 \left(\xi _1\right){}^2 b_0^4k^4\delta _1(k)b_0k^2\delta _2(k) \left(\xi _2\right){}^6+6 \left(\xi _1\right){}^2 b_0^4k^4\delta _2(k)b_0k^2\delta _1(k) \left(\xi _2\right){}^6
		\\ 
		&+24 \left(\xi _1\right){}^2 b_0^4k^5\delta _1(k)b_0k\delta _2(k) \left(\xi _2\right){}^6
		+24 \left(\xi _1\right){}^2 b_0^4k^5\delta _2(k)b_0k\delta _1(k) \left(\xi _2\right){}^6-2 b_0^3k^4\delta _1(k)\delta _2(k) \left(\xi _2\right){}^4
		\\
		&-3 b_0^3k^4\delta _2(k)\delta _1(k) \left(\xi _2\right){}^4
		-6 \left(\xi _1\right){}^2 b_0^4k^4\delta _1(k)\delta _2(k) \left(\xi _2\right){}^4
		-6 \left(\xi _1\right){}^2 b_0^4k^4\delta _2(k)\delta _1(k) \left(\xi _2\right){}^4
		\\
		&+b_0\delta _1(k)b_0^2k^4\delta _2(k) \left(\xi _2\right){}^4
		-\frac{1}{2} b_0\delta _2(k)b_0^2k^4\delta _1(k) \left(\xi _2\right){}^4-5 b_0^2k^4\delta _1(k)b_0\delta _2(k) \left(\xi _2\right){}^4
		\\
		&-6 b_0^2k^4\delta _2(k)b_0\delta _1(k) \left(\xi _2\right){}^4
		-5 \left(\xi _1\right){}^2 b_0^2k^4\delta _2(k)b_0^2\delta _1(k) \left(\xi _2\right){}^4-21 \left(\xi _1\right){}^2 b_0^3k^4\delta _1(k)b_0\delta _2(k) \left(\xi _2\right){}^4
		\\
		&-16 \left(\xi _1\right){}^2 b_0^3k^4\delta _2(k)b_0\delta _1(k) \left(\xi _2\right){}^4+b_0k\delta _1(k)b_0^2k^3\delta _2(k) \left(\xi _2\right){}^4-2 b_0k\delta _2(k)b_0^2k^3\delta _1(k) \left(\xi _2\right){}^4
		\\
		&-\frac{7}{2} b_0k^2\delta _2(k)b_0^2k^2\delta _1(k) \left(\xi _2\right){}^4-2 b_0k^3\delta _2(k)b_0^2k\delta _1(k) \left(\xi _2\right){}^4-2 \left(\xi _1\right){}^2 b_0^2k^2\delta _1(k)b_0^2k^2\delta _2(k) \left(\xi _2\right){}^4
		\\
		&-4 \left(\xi _1\right){}^2 b_0^2k^2\delta _2(k)b_0^2k^2\delta _1(k) \left(\xi _2\right){}^4-4 b_0^2k^3\delta _1(k)b_0k\delta _2(k) \left(\xi _2\right){}^4-2 \left(\xi _1\right){}^2 b_0^2k^3\delta _1(k)b_0^2k\delta _2(k) \left(\xi _2\right){}^4
		\\
		&-5 b_0^2k^3\delta _2(k)b_0k\delta _1(k) \left(\xi _2\right){}^4-9 \left(\xi _1\right){}^2 b_0^2k^3\delta _2(k)b_0^2k\delta _1(k) \left(\xi _2\right){}^4-4 \left(\xi _1\right){}^2 b_0^3k^2\delta _1(k)b_0k^2\delta _2(k) \left(\xi _2\right){}^4
		\\
		&-3 \left(\xi _1\right){}^2 b_0^3k^2\delta _2(k)b_0k^2\delta _1(k) \left(\xi _2\right){}^4-19 \left(\xi _1\right){}^2 b_0^3k^3\delta _1(k)b_0k\delta _2(k) \left(\xi _2\right){}^4
		\\
		&-13 \left(\xi _1\right){}^2 b_0^3k^3\delta _2(k)b_0k\delta _1(k) \left(\xi _2\right){}^4+2 \left(\xi _1\right){}^2 b_0^3k^2\delta _1(k)\delta _2(k) \left(\xi _2\right){}^2+3 \left(\xi _1\right){}^2 b_0^3k^2\delta _2(k)\delta _1(k) \left(\xi _2\right){}^2
		\\
		&+2 b_0k^2\delta _2(k)b_0\delta _1(k) \left(\xi _2\right){}^2+\frac{3}{2} \left(\xi _1\right){}^2 b_0k^2\delta _2(k)b_0^2\delta _1(k) \left(\xi _2\right){}^2-b_0\delta _1(k)b_0k^2\delta _2(k) \left(\xi _2\right){}^2
		\\
		&-\left(\xi _1\right){}^2 b_0\delta _1(k)b_0^2k^2\delta _2(k) \left(\xi _2\right){}^2+\frac{1}{2} b_0\delta _2(k)b_0k^2\delta _1(k) \left(\xi _2\right){}^2+\frac{1}{2} \left(\xi _1\right){}^2 b_0\delta _2(k)b_0^2k^2\delta _1(k) \left(\xi _2\right){}^2
		\\
		&+3 \left(\xi _1\right){}^2 b_0^2k^2\delta _1(k)b_0\delta _2(k) \left(\xi _2\right){}^2+\left(\xi _1\right){}^2 b_0^2k^2\delta _2(k)b_0\delta _1(k) \left(\xi _2\right){}^2-\frac{1}{2} b_0k\delta _1(k)b_0k\delta _2(k) \left(\xi _2\right){}^2
		\\
		&-\left(\xi _1\right){}^2 b_0k\delta _1(k)b_0^2k\delta _2(k) \left(\xi _2\right){}^2+\frac{3}{2} b_0k\delta _2(k)b_0k\delta _1(k) \left(\xi _2\right){}^2+2 \left(\xi _1\right){}^2 b_0k\delta _2(k)b_0^2k\delta _1(k) \left(\xi _2\right){}^2
		\\
		&+\frac{1}{2} \left(\xi _1\right){}^2 b_0\delta _1(k)b_0\delta _2(k)\Big) \\
	+ u_2 v_1
		&\Big(\left(-24 \left(\xi _2\right){}^4 \left(\xi _1\right){}^2 b_0^4k^5+4 \left(\xi _2\right){}^4 b_0^3k^5+8 \left(\xi _2\right){}^2 \left(\xi _1\right){}^4 b_0^4k^3+7 \left(\xi _2\right){}^2 \left(\xi _1\right){}^2 b_0^3k^3-3 \left(\xi _2\right){}^2 b_0^2k^3 \right.
		\\
		&\left. -\left(\xi _1\right){}^4 b_0^3k+\left(\xi _1\right){}^2 b_0^2k\right)\delta _1\left(\delta _2(k)\right) + \Big(-2 b_0^3k^6\delta _1(k)b_0\delta _2(k) \left(\xi _2\right){}^6-4 b_0^3k^6\delta _2(k)b_0\delta _1(k) \left(\xi _2\right){}^6
		\\
		& +2 \left(\xi _1\right){}^2 b_0^3k^6\delta _2(k)b_0^2\delta _1(k) \left(\xi _2\right){}^6 +18 \left(\xi _1\right){}^2 b_0^4k^6\delta _1(k)b_0\delta _2(k) \left(\xi _2\right){}^6+18 \left(\xi _1\right){}^2 b_0^4k^6\delta _2(k)b_0\delta _1(k) \left(\xi _2\right){}^6
		\\
		&-2 b_0^2k^2\delta _1(k)b_0^2k^4\delta _2(k) \left(\xi _2\right){}^6-4 b_0^2k^3\delta _1(k)b_0^2k^3\delta _2(k) \left(\xi _2\right){}^6-2 b_0^2k^4\delta _1(k)b_0^2k^2\delta _2(k) \left(\xi _2\right){}^6
		\\
		&+4 \left(\xi _1\right){}^2 b_0^3k^2\delta _1(k)b_0^2k^4\delta _2(k) \left(\xi _2\right){}^6+10 \left(\xi _1\right){}^2 b_0^3k^3\delta _1(k)b_0^2k^3\delta _2(k) \left(\xi _2\right){}^6
		\\
		& +4 \left(\xi _1\right){}^2 b_0^3k^3\delta _2(k)b_0^2k^3\delta _1(k) \left(\xi _2\right){}^6-2 b_0^3k^4\delta _1(k)b_0k^2\delta _2(k) \left(\xi _2\right){}^6
		+8 \left(\xi _1\right){}^2 b_0^3k^4\delta _1(k)b_0^2k^2\delta _2(k) \left(\xi _2\right){}^6
		\\
		&+10 \left(\xi _1\right){}^2 b_0^3k^4\delta _2(k)b_0^2k^2\delta _1(k) \left(\xi _2\right){}^6-4 b_0^3k^5\delta _1(k)b_0k\delta _2(k) \left(\xi _2\right){}^6+2 \left(\xi _1\right){}^2 b_0^3k^5\delta _1(k)b_0^2k\delta _2(k) \left(\xi _2\right){}^6
		\\
		&-4 b_0^3k^5\delta _2(k)b_0k\delta _1(k) \left(\xi _2\right){}^6+8 \left(\xi _1\right){}^2 b_0^3k^5\delta _2(k)b_0^2k\delta _1(k) \left(\xi _2\right){}^6+6 \left(\xi _1\right){}^2 b_0^4k^4\delta _1(k)b_0k^2\delta _2(k) \left(\xi _2\right){}^6
		\\
		&+6 \left(\xi _1\right){}^2 b_0^4k^4\delta _2(k)b_0k^2\delta _1(k) \left(\xi _2\right){}^6+24 \left(\xi _1\right){}^2 b_0^4k^5\delta _1(k)b_0k\delta _2(k) \left(\xi _2\right){}^6
		\\
		&+24 \left(\xi _1\right){}^2 b_0^4k^5\delta _2(k)b_0k\delta _1(k) \left(\xi _2\right){}^6+2 b_0^3k^4\delta _1(k)\delta _2(k) \left(\xi _2\right){}^4-6 \left(\xi _1\right){}^2 b_0^4k^4\delta _1(k)\delta _2(k) \left(\xi _2\right){}^4
		\\
		&-6 \left(\xi _1\right){}^2 b_0^4k^4\delta _2(k)\delta _1(k) \left(\xi _2\right){}^4+3 b_0^2k^4\delta _1(k)b_0\delta _2(k) \left(\xi _2\right){}^4+4 b_0^2k^4\delta _2(k)b_0\delta _1(k) \left(\xi _2\right){}^4
		\\
		&-2 \left(\xi _1\right){}^2 b_0^2k^4\delta _2(k)b_0^2\delta _1(k) \left(\xi _2\right){}^4-17 \left(\xi _1\right){}^2 b_0^3k^4\delta _1(k)b_0\delta _2(k) \left(\xi _2\right){}^4-15 \left(\xi _1\right){}^2 b_0^3k^4\delta _2(k)b_0\delta _1(k) \left(\xi _2\right){}^4
		\\
		&-2 \left(\xi _1\right){}^4 b_0^3k^4\delta _2(k)b_0^2\delta _1(k) \left(\xi _2\right){}^4-6 \left(\xi _1\right){}^4 b_0^4k^4\delta _1(k)b_0\delta _2(k) \left(\xi _2\right){}^4-6 \left(\xi _1\right){}^4 b_0^4k^4\delta _2(k)b_0\delta _1(k) \left(\xi _2\right){}^4
		\\
		&+2 b_0k\delta _1(k)b_0^2k^3\delta _2(k) \left(\xi _2\right){}^4+2 b_0k^2\delta _1(k)b_0^2k^2\delta _2(k) \left(\xi _2\right){}^4-\left(\xi _1\right){}^2 b_0^2k\delta _2(k)b_0^2k^3\delta _1(k) \left(\xi _2\right){}^4
		\\
		&+4 b_0^2k^2\delta _1(k)b_0k^2\delta _2(k) \left(\xi _2\right){}^4+b_0^2k^2\delta _2(k)b_0k^2\delta _1(k) \left(\xi _2\right){}^4-5 \left(\xi _1\right){}^2 b_0^2k^2\delta _2(k)b_0^2k^2\delta _1(k) \left(\xi _2\right){}^4
		\\
		&+7 b_0^2k^3\delta _1(k)b_0k\delta _2(k) \left(\xi _2\right){}^4+3 b_0^2k^3\delta _2(k)b_0k\delta _1(k) \left(\xi _2\right){}^4-6 \left(\xi _1\right){}^2 b_0^2k^3\delta _2(k)b_0^2k\delta _1(k) \left(\xi _2\right){}^4
		\\
		&-2 \left(\xi _1\right){}^4 b_0^3k\delta _1(k)b_0^2k^3\delta _2(k) \left(\xi _2\right){}^4-6 \left(\xi _1\right){}^2 b_0^3k^2\delta _1(k)b_0k^2\delta _2(k) \left(\xi _2\right){}^4
		\\
		&-4 \left(\xi _1\right){}^4 b_0^3k^2\delta _1(k)b_0^2k^2\delta _2(k) \left(\xi _2\right){}^4-2 \left(\xi _1\right){}^2 b_0^3k^2\delta _2(k)b_0k^2\delta _1(k) \left(\xi _2\right){}^4
		\\
		&-2 \left(\xi _1\right){}^4 b_0^3k^2\delta _2(k)b_0^2k^2\delta _1(k) \left(\xi _2\right){}^4-17 \left(\xi _1\right){}^2 b_0^3k^3\delta _1(k)b_0k\delta _2(k) \left(\xi _2\right){}^4
		\\
		&-2 \left(\xi _1\right){}^4 b_0^3k^3\delta _1(k)b_0^2k\delta _2(k) \left(\xi _2\right){}^4-11 \left(\xi _1\right){}^2 b_0^3k^3\delta _2(k)b_0k\delta _1(k) \left(\xi _2\right){}^4
		\\
		&-4 \left(\xi _1\right){}^4 b_0^3k^3\delta _2(k)b_0^2k\delta _1(k) \left(\xi _2\right){}^4-2 \left(\xi _1\right){}^4 b_0^4k^2\delta _1(k)b_0k^2\delta _2(k) \left(\xi _2\right){}^4
		\\
		&-2 \left(\xi _1\right){}^4 b_0^4k^2\delta _2(k)b_0k^2\delta _1(k) \left(\xi _2\right){}^4-8 \left(\xi _1\right){}^4 b_0^4k^3\delta _1(k)b_0k\delta _2(k) \left(\xi _2\right){}^4
		\\
		&-8 \left(\xi _1\right){}^4 b_0^4k^3\delta _2(k)b_0k\delta _1(k) \left(\xi _2\right){}^4-2 b_0^2k^2\delta _1(k)\delta _2(k) \left(\xi _2\right){}^2-b_0^2k^2\delta _2(k)\delta _1(k) \left(\xi _2\right){}^2
		\\
		&+2 \left(\xi _1\right){}^2 b_0^3k^2\delta _1(k)\delta _2(k) \left(\xi _2\right){}^2+2 \left(\xi _1\right){}^2 b_0^3k^2\delta _2(k)\delta _1(k) \left(\xi _2\right){}^2+2 \left(\xi _1\right){}^4 b_0^4k^2\delta _1(k)\delta _2(k) \left(\xi _2\right){}^2
		\\
		&+2 \left(\xi _1\right){}^4 b_0^4k^2\delta _2(k)\delta _1(k) \left(\xi _2\right){}^2-b_0k^2\delta _1(k)b_0\delta _2(k) \left(\xi _2\right){}^2-\left(\xi _1\right){}^2 b_0^2k^2\delta _1(k)b_0\delta _2(k) \left(\xi _2\right){}^2
		\\
		&+\left(\xi _1\right){}^4 b_0^2k^2\delta _2(k)b_0^2\delta _1(k) \left(\xi _2\right){}^2+5 \left(\xi _1\right){}^4 b_0^3k^2\delta _1(k)b_0\delta _2(k) \left(\xi _2\right){}^2+3 \left(\xi _1\right){}^4 b_0^3k^2\delta _2(k)b_0\delta _1(k) \left(\xi _2\right){}^2
		\\
		&-2 b_0k\delta _1(k)b_0k\delta _2(k) \left(\xi _2\right){}^2-\left(\xi _1\right){}^2 b_0^2k\delta _1(k)b_0k\delta _2(k) \left(\xi _2\right){}^2+\left(\xi _1\right){}^4 b_0^2k\delta _2(k)b_0^2k\delta _1(k) \left(\xi _2\right){}^2
		\\
		&+3 \left(\xi _1\right){}^4 b_0^3k\delta _1(k)b_0k\delta _2(k) \left(\xi _2\right){}^2+\left(\xi _1\right){}^4 b_0^3k\delta _2(k)b_0k\delta _1(k) \left(\xi _2\right){}^2\Big)\Big) \\
+ u_2v_2 \Big(
		&\left(4 \left(\xi _2\right){}^6 b_0^4k^7-2 \left(\xi _2\right){}^4 b_0^3k^5-12 \left(\xi _1\right){}^2 \left(\xi _2\right){}^4 b_0^4k^5+2 \left(\xi _1\right){}^2 \left(\xi _2\right){}^2 b_0^3k^3\right)\delta _1\left(\delta _1(k)\right)
		\\
		&\left(-4 \left(\xi _2\right){}^6 b_0^4k^9+7 \left(\xi _2\right){}^4 b_0^3k^7+12 \left(\xi _1\right){}^2 \left(\xi _2\right){}^4 b_0^4k^7 -3 \left(\xi _2\right){}^2 b_0^2k^5 - 9 \left(\xi _1\right){}^2 \left(\xi _2\right){}^2 b_0^3k^5 \right.
		\\
		&
		\left.+\left(\xi _1\right){}^2 b_0^2k^3\right)\delta _2\left(\delta _2(k)\right)
		+
		\Big(-6 b_0^4k^8\delta _1(k)b_0\delta _1(k) \left(\xi _2\right){}^8-2 b_0^3k^5\delta _1(k)b_0^2k^3\delta _1(k) \left(\xi _2\right){}^8
		\\
		&-4 b_0^3k^6\delta _1(k)b_0^2k^2\delta _1(k) \left(\xi _2\right){}^8-2 b_0^3k^7\delta _1(k)b_0^2k\delta _1(k) \left(\xi _2\right){}^8-2 b_0^4k^6\delta _1(k)b_0k^2\delta _1(k) \left(\xi _2\right){}^8
		\\
		&-8 b_0^4k^7\delta _1(k)b_0k\delta _1(k) \left(\xi _2\right){}^8+2 b_0^4k^6\delta _1(k){}^2 \left(\xi _2\right){}^6+8 b_0^3k^6\delta _1(k)b_0\delta _1(k) \left(\xi _2\right){}^6
		\\
		&+4 \left(\xi _1\right){}^2 b_0^3k^6\delta _1(k)b_0^2\delta _1(k) \left(\xi _2\right){}^6+18 \left(\xi _1\right){}^2 b_0^4k^6\delta _1(k)b_0\delta _1(k) \left(\xi _2\right){}^6+2 b_0^2k^4\delta _1(k)b_0^2k^2\delta _1(k) \left(\xi _2\right){}^6
		\\
		&+2 b_0^2k^5\delta _1(k)b_0^2k\delta _1(k) \left(\xi _2\right){}^6+2 \left(\xi _1\right){}^2 b_0^3k^3\delta _1(k)b_0^2k^3\delta _1(k) \left(\xi _2\right){}^6+8 \left(\xi _1\right){}^2 b_0^3k^4\delta _1(k)b_0^2k^2\delta _1(k) \left(\xi _2\right){}^6
		\\
		&+6 b_0^3k^5\delta _1(k)b_0k\delta _1(k) \left(\xi _2\right){}^6+10 \left(\xi _1\right){}^2 b_0^3k^5\delta _1(k)b_0^2k\delta _1(k) \left(\xi _2\right){}^6+6 \left(\xi _1\right){}^2 b_0^4k^4\delta _1(k)b_0k^2\delta _1(k) \left(\xi _2\right){}^6
		\\
		&+24 \left(\xi _1\right){}^2 b_0^4k^5\delta _1(k)b_0k\delta _1(k) \left(\xi _2\right){}^6-6 \left(\xi _1\right){}^2 b_0^4k^4\delta _1(k){}^2 \left(\xi _2\right){}^4-2 b_0^2k^4\delta _1(k)b_0\delta _1(k) \left(\xi _2\right){}^4
		\\
		&-2 \left(\xi _1\right){}^2 b_0^2k^4\delta _1(k)b_0^2\delta _1(k) \left(\xi _2\right){}^4-12 \left(\xi _1\right){}^2 b_0^3k^4\delta _1(k)b_0\delta _1(k) \left(\xi _2\right){}^4+b_0k\delta _1(k)b_0^2k^3\delta _1(k) \left(\xi _2\right){}^4
		\\
		&+b_0k^2\delta _1(k)b_0^2k^2\delta _1(k) \left(\xi _2\right){}^4-2 \left(\xi _1\right){}^2 b_0^2k^3\delta _1(k)b_0^2k\delta _1(k) \left(\xi _2\right){}^4-6 \left(\xi _1\right){}^2 b_0^3k^3\delta _1(k)b_0k\delta _1(k) \left(\xi _2\right){}^4
		\\
		&-\left(\xi _1\right){}^2 b_0k^2\delta _1(k)b_0^2\delta _1(k) \left(\xi _2\right){}^2-b_0k\delta _1(k)b_0k\delta _1(k) \left(\xi _2\right){}^2-\left(\xi _1\right){}^2 b_0k\delta _1(k)b_0^2k\delta _1(k) \left(\xi _2\right){}^2\Big)
		\\
		& +\Big(6 b_0^4k^{10}\delta _2(k)b_0\delta _2(k) \left(\xi _2\right){}^8+2 b_0^3k^6\delta _2(k)b_0^2k^4\delta _2(k) \left(\xi _2\right){}^8+4 b_0^3k^7\delta _2(k)b_0^2k^3\delta _2(k) \left(\xi _2\right){}^8
		\\
		&+2 b_0^3k^8\delta _2(k)b_0^2k^2\delta _2(k) \left(\xi _2\right){}^8+2 b_0^4k^8\delta _2(k)b_0k^2\delta _2(k) \left(\xi _2\right){}^8+8 b_0^4k^9\delta _2(k)b_0k\delta _2(k) \left(\xi _2\right){}^8
		\\
		&-2 b_0^4k^8\delta _2(k){}^2 \left(\xi _2\right){}^6-13 b_0^3k^8\delta _2(k)b_0\delta _2(k) \left(\xi _2\right){}^6-18 \left(\xi _1\right){}^2 b_0^4k^8\delta _2(k)b_0\delta _2(k) \left(\xi _2\right){}^6
		\\
		&-3 b_0^2k^4\delta _2(k)b_0^2k^4\delta _2(k) \left(\xi _2\right){}^6-7 b_0^2k^5\delta _2(k)b_0^2k^3\delta _2(k) \left(\xi _2\right){}^6-4 b_0^2k^6\delta _2(k)b_0^2k^2\delta _2(k) \left(\xi _2\right){}^6 
		\\
		&-2 \left(\xi _1\right){}^2 b_0^3k^4\delta _2(k)b_0^2k^4\delta _2(k) \left(\xi _2\right){}^6-8 \left(\xi _1\right){}^2 b_0^3k^5\delta _2(k)b_0^2k^3\delta _2(k) \left(\xi _2\right){}^6-7 b_0^3k^6\delta _2(k)b_0k^2\delta _2(k) \left(\xi _2\right){}^6
		\\
		&-10 \left(\xi _1\right){}^2 b_0^3k^6\delta _2(k)b_0^2k^2\delta _2(k) \left(\xi _2\right){}^6-18 b_0^3k^7\delta _2(k)b_0k\delta _2(k) \left(\xi _2\right){}^6-4 \left(\xi _1\right){}^2 b_0^3k^7\delta _2(k)b_0^2k\delta _2(k) \left(\xi _2\right){}^6
		\\
		&-6 \left(\xi _1\right){}^2 b_0^4k^6\delta _2(k)b_0k^2\delta _2(k) \left(\xi _2\right){}^6-24 \left(\xi _1\right){}^2 b_0^4k^7\delta _2(k)b_0k\delta _2(k) \left(\xi _2\right){}^6+5 b_0^3k^6\delta _2(k){}^2 \left(\xi _2\right){}^4
		\\
		&+6 \left(\xi _1\right){}^2 b_0^4k^6\delta _2(k){}^2 \left(\xi _2\right){}^4+8 b_0^2k^6\delta _2(k)b_0\delta _2(k) \left(\xi _2\right){}^4+25 \left(\xi _1\right){}^2 b_0^3k^6\delta _2(k)b_0\delta _2(k) \left(\xi _2\right){}^4
		\\
		&+\frac{1}{2} b_0k\delta _2(k)b_0^2k^5\delta _2(k) \left(\xi _2\right){}^4+2 b_0k^2\delta _2(k)b_0^2k^4\delta _2(k) \left(\xi _2\right){}^4+\frac{7}{2} b_0k^3\delta _2(k)b_0^2k^3\delta _2(k) \left(\xi _2\right){}^4
		\\
		&+2 b_0k^4\delta _2(k)b_0^2k^2\delta _2(k) \left(\xi _2\right){}^4+4 \left(\xi _1\right){}^2 b_0^2k^3\delta _2(k)b_0^2k^3\delta _2(k) \left(\xi _2\right){}^4+6 b_0^2k^4\delta _2(k)b_0k^2\delta _2(k) \left(\xi _2\right){}^4
		\\
		&+9 \left(\xi _1\right){}^2 b_0^2k^4\delta _2(k)b_0^2k^2\delta _2(k) \left(\xi _2\right){}^4+\frac{25}{2} b_0^2k^5\delta _2(k)b_0k\delta _2(k) \left(\xi _2\right){}^4+5 \left(\xi _1\right){}^2 b_0^2k^5\delta _2(k)b_0^2k\delta _2(k) \left(\xi _2\right){}^4
		\\
		&+7 \left(\xi _1\right){}^2 b_0^3k^4\delta _2(k)b_0k^2\delta _2(k) \left(\xi _2\right){}^4+26 \left(\xi _1\right){}^2 b_0^3k^5\delta _2(k)b_0k\delta _2(k) \left(\xi _2\right){}^4-3 b_0^2k^4\delta _2(k){}^2 \left(\xi _2\right){}^2
		\\
		&-5 \left(\xi _1\right){}^2 b_0^3k^4\delta _2(k){}^2 \left(\xi _2\right){}^2-b_0k^4\delta _2(k)b_0\delta _2(k) \left(\xi _2\right){}^2-8 \left(\xi _1\right){}^2 b_0^2k^4\delta _2(k)b_0\delta _2(k) \left(\xi _2\right){}^2
		\\
		&-\frac{1}{2} b_0k\delta _2(k)b_0k^3\delta _2(k) \left(\xi _2\right){}^2
		-\frac{1}{2} \left(\xi _1\right){}^2 b_0k\delta _2(k)b_0^2k^3\delta _2(k) \left(\xi _2\right){}^2-\frac{7}{4} b_0k^2\delta _2(k)b_0k^2\delta _2(k) \left(\xi _2\right){}^2
		\\
		&-2 \left(\xi _1\right){}^2 b_0k^2\delta _2(k)b_0^2k^2\delta _2(k) \left(\xi _2\right){}^2
		-\frac{11}{4} b_0k^3\delta _2(k)b_0k\delta _2(k) \left(\xi _2\right){}^2-\frac{3}{2} \left(\xi _1\right){}^2 b_0k^3\delta _2(k)b_0^2k\delta _2(k) \left(\xi _2\right){}^2
		\\
		&-\frac{11}{2} \left(\xi _1\right){}^2 b_0^2k^3\delta _2(k)b_0k\delta _2(k) \left(\xi _2\right){}^2
		+\frac{3}{4} \left(\xi _1\right){}^2 b_0k^2\delta _2(k)b_0\delta _2(k)+\frac{1}{4} \left(\xi _1\right){}^2 b_0k\delta _2(k)b_0k\delta _2(k)\Big)\Big){\Bigg)}
\end{align*}

\section{Appendix: the integrals}

In the next step we compute all integrals for all terms that appear in the symbol above. To make the procudure clear, we split the
computations in several steps and tabularize the results. The results were verified using Wolfram's Mathematica\footnote{To code is available from the authors}. 

The sum of the evaluated integrals for each table vanishes. 

\begin{longtable}{|c|c|c|c|}
	\caption{Terms at $u_1v_1$ of the form $C \,\xi_1^{2\alpha} \xi_2^{2\beta}b_0^m k^a \delta_1(k)b_0^n k^b \delta_1(k)$}\label{tab1}\\
	\hline
	\multicolumn{1}{|c|}{No.} &
	\multicolumn{1}{|c|}{Terms} &
	\multicolumn{1}{|c|}{\parbox{5cm}{\centering $\text{C} \times F(s, m, n,\beta, a, b)$}} &
	\multicolumn{1}{|c|}{\parbox{5cm}{\centering Result}} \\
	\hline
	\endfirsthead
	\multicolumn{4}{|c|}{{\bfseries \tablename\ \thetable{} -- continued from previous page}} \\
	\hline
	\multicolumn{1}{|c|}{SL No.} &
	\multicolumn{1}{|c|}{Terms} &
	\multicolumn{1}{|c|}{\parbox{5cm}{\centering $\text{Coeff} \times F(s, m, n,\beta, a, b)$}} &
	\multicolumn{1}{|c|}{\parbox{5cm}{\centering Result}} \\
	\hline\endhead
	\hline
	\endfoot
	\hline
	\endlastfoot
	1 & $2\, \xi _1{}^2 \xi _2{}^2 b_0{}^3\delta _1(k){}^2$ & \parbox{5cm}{\centering $2\, F(s,3,0,1,0,0)$} & \parbox{5cm}{\centering $\frac{\pi }{8 k^3}$} \\ \hline
	2 & $-2\, \xi _1{}^4 \xi _2{}^2 b_0{}^4\delta _1(k){}^2$ & \parbox{5cm}{\centering $-2\, F(s,4,0,1,0,0)$} & \parbox{5cm}{\centering $-\frac{\pi }{16 k^3}$} \\ \hline
	3 & $-2\, \xi _2{}^4 b_0{}^3k^2\delta _1(k){}^2$ & \parbox{5cm}{\centering $-2\, F(s,3,0,2,2,0)$} & \parbox{5cm}{\centering $-\frac{3 \pi }{8 k^3}$} \\ \hline
	4 & $6\, \xi _1{}^2 \xi _2{}^4 b_0{}^4k^2\delta _1(k){}^2$ & \parbox{5cm}{\centering $6\, F(s,4,0,2,2,0)$} & \parbox{5cm}{\centering $\frac{3 \pi }{16 k^3}$} \\ \hline
	5 & $2\, \xi _2{}^2 b_0\delta _1(k)b_0\delta _1(k)$ & \parbox{5cm}{\centering $2\, F(s,1,1,1,0,0)$} & \parbox{5cm}{\centering $\frac{2 \pi }{k^3 \left(2 s^2+2 s\right)}$} \\ \hline
	6 & $-2\, \xi _1{}^4 \xi _2{}^2 b_0{}^3\delta _1(k)b_0\delta _1(k)$ & \parbox{5cm}{\centering $-2\, F(s,3,1,1,0,0)$} & \parbox{5cm}{\centering $-\frac{\pi  (3 s+1)}{8 k^3 (s+1)^3}$} \\ \hline
	7 & $-2\, \xi _2{}^4 b_0\delta _1(k)b_0{}^2k^2\delta _1(k)$ & \parbox{5cm}{\centering $-2\, F(s,1,2,2,0,2)$} & \parbox{5cm}{\centering $-\frac{\pi  (2 s+1)}{2 k^3 s (s+1)^2}$} \\ \hline
	8 & $-8\, \xi _2{}^4 b_0{}^2k^2\delta _1(k)b_0\delta _1(k)$ & \parbox{5cm}{\centering $-8\, F(s,2,1,2,2,0)$} & \parbox{5cm}{\centering $-\frac{2 \pi  (s+2)}{k^3 s (s+1)^2}$} \\
\end{longtable}

\begin{longtable}{|c|c|c|c|}
	\caption{Terms at $u_1v_1$ of the form $C \,\xi_1^{2\alpha} \xi_2^{2\beta}b_0^m k^a \delta_2(k)b_0^n k^b \delta_2(k)$}\label{tab2}\\
	\hline
 \multicolumn{1}{|c|}{No.} &
	\multicolumn{1}{|c|}{Terms} &
	\multicolumn{1}{|c|}{\parbox{5cm}{\centering $\text{C} \times F(s, m, n,\beta, a, b)$}} &
	\multicolumn{1}{|c|}{\parbox{5cm}{\centering Result}} \\
	\hline
	\endfirsthead
	\multicolumn{4}{|c|}{{\bfseries \tablename\ \thetable{} -- continued from previous page}} \\
	\hline
	\multicolumn{1}{|c|}{SL No.} &
	\multicolumn{1}{|c|}{Terms} &
	\multicolumn{1}{|c|}{\parbox{5cm}{\centering $\text{Coeff} \times F(s, m, n,\beta, a, b)$}} &
	\multicolumn{1}{|c|}{\parbox{5cm}{\centering Result}} \\
	\hline\endhead
	\hline
	\endfoot
	\hline
	\endlastfoot
	1 & $6\, \xi _1{}^2 \xi _2{}^2 b_0{}^3k^2\delta _2(k){}^2$ & \parbox{5cm}{\centering $6\, F(s,3,0,1,2,0)$} & \parbox{5cm}{\centering $\frac{3 \pi }{8 k}$} \\ \hline
	2 & $2\, \xi _1{}^4 \xi _2{}^2 b_0{}^4k^2\delta _2(k){}^2$ & \parbox{5cm}{\centering $2\, F(s,4,0,1,2,0)$} & \parbox{5cm}{\centering $\frac{\pi }{16 k}$} \\ \hline
	3 & $-6\, \xi _1{}^2 \xi _2{}^4 b_0{}^4k^4\delta _2(k){}^2$ & \parbox{5cm}{\centering $-6\, F(s,4,0,2,4,0)$} & \parbox{5cm}{\centering $-\frac{3 \pi }{16 k}$} \\ \hline
	4 & $-\frac{1}{2} \xi _1{}^4 b_0{}^2\delta _2(k)b_0\delta _2(k)$ & \parbox{5cm}{\centering $-\frac{1}{2} F(s,2,1,0,0,0)$} & \parbox{5cm}{\centering $-\frac{\pi  (2 s+1)}{8 k (s+1)^2}$} \\ \hline
	5 & $\xi _1{}^2 \xi _2{}^2 b_0{}^2\delta _2(k)b_0k^2\delta _2(k)$ & \parbox{5cm}{\centering $F(s,2,1,1,0,2)$} & \parbox{5cm}{\centering $\frac{\pi  s^2}{4 k (s+1)^2}$} \\ \hline
	6 & $\xi _1{}^4 \xi _2{}^2 b_0{}^2\delta _2(k)b_0{}^2k^2\delta _2(k)$ & \parbox{5cm}{\centering $F(s,2,2,1,0,2)$} & \parbox{5cm}{\centering $\frac{\pi  s^2}{4 k (s+1)^3}$} \\ \hline
	7 & $-\xi _1{}^2 \xi _2{}^4 b_0{}^2\delta _2(k)b_0{}^2k^4\delta _2(k)$ & \parbox{5cm}{\centering $-F(s,2,2,2,0,4)$} & \parbox{5cm}{\centering $-\frac{\pi  s^3}{4 k (s+1)^3}$} \\ \hline
	8 & $\frac{5}{2} \xi _1{}^2 \xi _2{}^2 b_0{}^2k^2\delta _2(k)b_0\delta _2(k)$ & \parbox{5cm}{\centering $\frac{5}{2} F(s,2,1,1,2,0)$} & \parbox{5cm}{\centering $\frac{5 \pi }{8 k (s+1)^2}$} \\
\end{longtable}

\begin{longtable}{|c|c|c|c|}
	\caption{Terms at $u_1v_2$ of the form $C \,\xi_1^{2\alpha} \xi_2^{2\beta}b_0^m k^a \delta_1(k)b_0^n k^b \delta_2(k)$}\label{tab3}\\
	\hline
	\multicolumn{1}{|c|}{No.} &
\multicolumn{1}{|c|}{Terms} &
\multicolumn{1}{|c|}{\parbox{5cm}{\centering $\text{C} \times F(s, m, n,\beta, a, b)$}} &
\multicolumn{1}{|c|}{\parbox{5cm}{\centering Result}} \\
	\hline
	\endfirsthead
	\multicolumn{4}{|c|}{{\bfseries \tablename\ \thetable{} -- continued from previous page}} \\
	\hline
	\multicolumn{1}{|c|}{SL No.} &
	\multicolumn{1}{|c|}{Terms} &
	\multicolumn{1}{|c|}{\parbox{5cm}{\centering $\text{Coeff} \times F(s, m, n,\beta, a, b)$}} &
	\multicolumn{1}{|c|}{\parbox{5cm}{\centering Result}} \\
	\hline\endhead
	\hline
	\endfoot
	\hline
	\endlastfoot
	1 & $\frac{1}{2} \xi _1{}^2 b_0\delta _1(k)b_0\delta _2(k)$ & \parbox{5cm}{\centering $\frac{1}{2} F(s,1,1,0,0,0)$} & \parbox{5cm}{\centering $\frac{\pi }{2 k (2 s+2)}$} \\ \hline
	2 & $2\, \xi _1{}^2 \xi _2{}^2 b_0{}^3k^2\delta _1(k)\delta _2(k)$ & \parbox{5cm}{\centering $2\, F(s,3,0,1,2,0)$} & \parbox{5cm}{\centering $\frac{\pi }{8 k}$} \\ \hline
	3 & $3\, \xi _1{}^2 \xi _2{}^2 b_0{}^3k^2\delta _2(k)\delta _1(k)$ & \parbox{5cm}{\centering $3\, F(s,0,3,1,0,2)$} & \parbox{5cm}{\centering $\frac{3 \pi }{16 k s}$} \\ \hline
	4 & $-2\, \xi _2{}^4 b_0{}^3k^4\delta _1(k)\delta _2(k)$ & \parbox{5cm}{\centering $-2\, F(s,3,0,2,4,0)$} & \parbox{5cm}{\centering $-\frac{3 \pi }{8 k}$} \\ \hline
	5 & $-3\, \xi _2{}^4 b_0{}^3k^4\delta _2(k)\delta _1(k)$ & \parbox{5cm}{\centering $-3\, F(s,0,3,2,0,4)$} & \parbox{5cm}{\centering $-\frac{9 \pi }{16 k s}$} \\ \hline
	6 & $-6\, \xi _1{}^2 \xi _2{}^4 b_0{}^4k^4\delta _1(k)\delta _2(k)$ & \parbox{5cm}{\centering $-6\, F(s,4,0,2,4,0)$} & \parbox{5cm}{\centering $-\frac{3 \pi }{16 k}$} \\ \hline
	7 & $-6\, \xi _1{}^2 \xi _2{}^4 b_0{}^4k^4\delta _2(k)\delta _1(k)$ & \parbox{5cm}{\centering $-6\, F(s,0,4,2,0,4)$} & \parbox{5cm}{\centering $-\frac{3 \pi }{16 k s}$} \\ \hline
	8 & $2\, \xi _2{}^6 b_0{}^4k^6\delta _1(k)\delta _2(k)$ & \parbox{5cm}{\centering $2\, F(s,4,0,3,6,0)$} & \parbox{5cm}{\centering $\frac{5 \pi }{16 k}$} \\
\end{longtable}

\begin{longtable}{|c|c|c|c|}
	\caption{Terms at $u_2v_1$ of the form $C \,\xi_1^{2\alpha} \xi_2^{2\beta}b_0^m k^a \delta_1(k)b_0^n k^b \delta_2(k)$}\label{tab4}\\
	\hline
	\multicolumn{1}{|c|}{No.} &
\multicolumn{1}{|c|}{Terms} &
\multicolumn{1}{|c|}{\parbox{5cm}{\centering $\text{C} \times F(s, m, n,\beta, a, b)$}} &
\multicolumn{1}{|c|}{\parbox{5cm}{\centering Result}} \\
	\hline
	\endfirsthead
	\multicolumn{4}{|c|}{{\bfseries \tablename\ \thetable{} -- continued from previous page}} \\
	\hline
	\multicolumn{1}{|c|}{SL No.} &
	\multicolumn{1}{|c|}{Terms} &
	\multicolumn{1}{|c|}{\parbox{5cm}{\centering $\text{Coeff} \times F(s, m, n,\beta, a, b)$}} &
	\multicolumn{1}{|c|}{\parbox{5cm}{\centering Result}} \\
	\hline\endhead
	\hline
	\endfoot
	\hline
	\endlastfoot
	1 & $-2\, \xi _2{}^2 b_0{}^2k^2\delta _1(k)\delta _2(k)$ & \parbox{5cm}{\centering $-2\, F(s,2,0,1,2,0)$} & \parbox{5cm}{\centering $-\frac{\pi }{2 k}$} \\ \hline
	2 & $-\xi _2{}^2 b_0{}^2k^2\delta _2(k)\delta _1(k)$ & \parbox{5cm}{\centering $-F(s,0,2,1,0,2)$} & \parbox{5cm}{\centering $-\frac{\pi }{4 k s}$} \\ \hline
	3 & $2\, \xi _1{}^2 \xi _2{}^2 b_0{}^3k^2\delta _1(k)\delta _2(k)$ & \parbox{5cm}{\centering $2\, F(s,3,0,1,2,0)$} & \parbox{5cm}{\centering $\frac{\pi }{8 k}$} \\ \hline
	4 & $2\, \xi _1{}^2 \xi _2{}^2 b_0{}^3k^2\delta _2(k)\delta _1(k)$ & \parbox{5cm}{\centering $2\, F(s,0,3,1,0,2)$} & \parbox{5cm}{\centering $\frac{\pi }{8 k s}$} \\ \hline
	5 & $2\, \xi _2{}^4 b_0{}^3k^4\delta _1(k)\delta _2(k)$ & \parbox{5cm}{\centering $2\, F(s,3,0,2,4,0)$} & \parbox{5cm}{\centering $\frac{3 \pi }{8 k}$} \\ \hline
	6 & $2\, \xi _1{}^4 \xi _2{}^2 b_0{}^4k^2\delta _1(k)\delta _2(k)$ & \parbox{5cm}{\centering $2\, F(s,4,0,1,2,0)$} & \parbox{5cm}{\centering $\frac{\pi }{16 k}$} \\ \hline
	7 & $2\, \xi _1{}^4 \xi _2{}^2 b_0{}^4k^2\delta _2(k)\delta _1(k)$ & \parbox{5cm}{\centering $2\, F(s,0,4,1,0,2)$} & \parbox{5cm}{\centering $\frac{\pi }{16 k s}$} \\ \hline
	8 & $-6\, \xi _1{}^2 \xi _2{}^4 b_0{}^4k^4\delta _1(k)\delta _2(k)$ & \parbox{5cm}{\centering $-6\, F(s,4,0,2,4,0)$} & \parbox{5cm}{\centering $-\frac{3 \pi }{16 k}$} \\
\end{longtable}

\begin{longtable}{|c|c|c|c|}
	\caption{Terms at $u_2 v_2$ of the form $C \,\xi_1^{2\alpha} \xi_2^{2\beta}b_0^m k^a \delta_1(k)b_0^n k^b \delta_1(k)$}\label{tab5}\\
	\hline
	\multicolumn{1}{|c|}{No.} &
\multicolumn{1}{|c|}{Terms} &
\multicolumn{1}{|c|}{\parbox{5cm}{\centering $\text{C} \times F(s, m, n,\beta, a, b)$}} &
\multicolumn{1}{|c|}{\parbox{5cm}{\centering Result}} \\
	\hline
	\endfirsthead
	\multicolumn{4}{|c|}{{\bfseries \tablename\ \thetable{} -- continued from previous page}} \\
	\hline
	\multicolumn{1}{|c|}{SL No.} &
	\multicolumn{1}{|c|}{Terms} &
	\multicolumn{1}{|c|}{\parbox{5cm}{\centering $\text{Coeff} \times F(s, m, n,\beta, A, B)$}} &
	\multicolumn{1}{|c|}{\parbox{5cm}{\centering Result}} \\
	\hline\endhead
	\hline
	\endfoot
	\hline
	\endlastfoot
	1 & $-6\, \xi _1{}^2 \xi _2{}^4 b_0{}^4k^4\delta _1(k){}^2$ & \parbox{5cm}{\centering $-6\, F(s,4,0,2,4,0)$} & \parbox{5cm}{\centering $-\frac{3 \pi }{16 k}$} \\ \hline
	2 & $2\, \xi _2{}^6 b_0{}^4k^6\delta _1(k){}^2$ & \parbox{5cm}{\centering $2\, F(s,4,0,3,6,0)$} & \parbox{5cm}{\centering $\frac{5 \pi }{16 k}$} \\ \hline
	3 & $-\xi _1{}^2 \xi _2{}^2 b_0k^2\delta _1(k)b_0{}^2\delta _1(k)$ & \parbox{5cm}{\centering $-F(s,1,2,1,2,0)$} & \parbox{5cm}{\centering $-\frac{\pi }{4 k s (s+1)^2}$} \\ \hline
	4 & $-2\, \xi _2{}^4 b_0{}^2k^4\delta _1(k)b_0\delta _1(k)$ & \parbox{5cm}{\centering $-2\, F(s,2,1,2,4,0)$} & \parbox{5cm}{\centering $-\frac{\pi  (s+2)}{2 k s (s+1)^2}$} \\ \hline
	5 & $-2\, \xi _1{}^2 \xi _2{}^4 b_0{}^2k^4\delta _1(k)b_0{}^2\delta _1(k)$ & \parbox{5cm}{\centering $-2\, F(s,2,2,2,4,0)$} & \parbox{5cm}{\centering $-\frac{\pi }{2 k s (s+1)^3}$} \\ \hline
	6 & $-12\, \xi _1{}^2 \xi _2{}^4 b_0{}^3k^4\delta _1(k)b_0\delta _1(k)$ & \parbox{5cm}{\centering $-12\, F(s,3,1,2,4,0)$} & \parbox{5cm}{\centering $-\frac{3 \pi  (s+3)}{4 k (s+1)^3}$} \\ \hline
	7 & $8\, \xi _2{}^6 b_0{}^3k^6\delta _1(k)b_0\delta _1(k)$ & \parbox{5cm}{\centering $8\, F(s,3,1,3,6,0)$} & \parbox{5cm}{\centering $\frac{\pi  (3 s (s+3)+8)}{2 k s (s+1)^3}$} \\ \hline
	8 & $4\, \xi _1{}^2 \xi _2{}^6 b_0{}^3k^6\delta _1(k)b_0{}^2\delta _1(k)$ & \parbox{5cm}{\centering $4\, F(s,3,2,3,6,0)$} & \parbox{5cm}{\centering $\frac{\pi  (s+4)}{4 k s (s+1)^4}$} \\
\end{longtable}

\begin{longtable}{|c|c|c|c|}
	\caption{Terms at $u_2 v_2$ of the form $C \,\xi_1^{2\alpha} \xi_2^{2\beta}b_0^m k^a \delta_2(k)b_0^n k^b \delta_2(k)$}\label{tab6}\\
	\hline
	\multicolumn{1}{|c|}{No.} &
  \multicolumn{1}{|c|}{Terms} &
  \multicolumn{1}{|c|}{\parbox{5cm}{\centering $\text{C} \times F(s, m, n,\beta, a, b)$}} &
  \multicolumn{1}{|c|}{\parbox{5cm}{\centering Result}} \\
	\hline
	\endfirsthead
	\multicolumn{4}{|c|}{{\bfseries \tablename\ \thetable{} -- continued from previous page}} \\
	\hline
	\multicolumn{1}{|c|}{SL No.} &
	\multicolumn{1}{|c|}{Terms} &
	\multicolumn{1}{|c|}{\parbox{5cm}{\centering $\text{Coeff} \times F(s, m, n,\beta, A, B)$}} &
	\multicolumn{1}{|c|}{\parbox{5cm}{\centering Result}} \\
	\hline\endhead
	\hline
	\endfoot
	\hline
	\endlastfoot
	1 & $-3\, \xi _2{}^2 b_0{}^2k^4\delta _2(k){}^2$ & \parbox{5cm}{\centering $-3\, F(s,2,0,1,4,0)$} & \parbox{5cm}{\centering $-\frac{3 \pi  k}{4}$} \\ \hline
	2 & $-5\, \xi _1{}^2 \xi _2{}^2 b_0{}^3k^4\delta _2(k){}^2$ & \parbox{5cm}{\centering $-5\, F(s,3,0,1,4,0)$} & \parbox{5cm}{\centering $-\frac{5 \pi  k}{16}$} \\ \hline
	3 & $5\, \xi _2{}^4 b_0{}^3k^6\delta _2(k){}^2$ & \parbox{5cm}{\centering $5\, F(s,3,0,2,6,0)$} & \parbox{5cm}{\centering $\frac{15 \pi  k}{16}$} \\ \hline
	4 & $6\, \xi _1{}^2 \xi _2{}^4 b_0{}^4k^6\delta _2(k){}^2$ & \parbox{5cm}{\centering $6\, F(s,4,0,2,6,0)$} & \parbox{5cm}{\centering $\frac{3 \pi  k}{16}$} \\ \hline
	5 & $-2\, \xi _2{}^6 b_0{}^4k^8\delta _2(k){}^2$ & \parbox{5cm}{\centering $-2\, F(s,4,0,3,8,0)$} & \parbox{5cm}{\centering $-\frac{5 \pi  k}{16}$} \\ \hline
	6 & $\frac{3}{4} \xi _1{}^2 b_0k^2\delta _2(k)b_0\delta _2(k)$ & \parbox{5cm}{\centering $\frac{3}{4} F(s,1,1,0,2,0)$} & \parbox{5cm}{\centering $\frac{3 \pi  k}{4 (2 s+2)}$} \\ \hline
	7 & $-\xi _2{}^2 b_0k^4\delta _2(k)b_0\delta _2(k)$ & \parbox{5cm}{\centering $-F(s,1,1,1,4,0)$} & \parbox{5cm}{\centering $-\frac{\pi  k}{2 s^2+2 s}$} \\ \hline
	8 & $-8\, \xi _1{}^2 \xi _2{}^2 b_0{}^2k^4\delta _2(k)b_0\delta _2(k)$ & \parbox{5cm}{\centering $-8\, F(s,2,1,1,4,0)$} & \parbox{5cm}{\centering $-\frac{2 \pi  k}{(s+1)^2}$} \\
\end{longtable}


\begin{thebibliography}{15}
\ifx \bisbn   \undefined \def \bisbn  #1{ISBN #1}\fi
\ifx \binits  \undefined \def \binits#1{#1}\fi
\ifx \bauthor  \undefined \def \bauthor#1{#1}\fi
\ifx \batitle  \undefined \def \batitle#1{#1}\fi
\ifx \bjtitle  \undefined \def \bjtitle#1{#1}\fi
\ifx \bvolume  \undefined \def \bvolume#1{\textbf{#1}}\fi
\ifx \byear  \undefined \def \byear#1{#1}\fi
\ifx \bissue  \undefined \def \bissue#1{#1}\fi
\ifx \bfpage  \undefined \def \bfpage#1{#1}\fi
\ifx \blpage  \undefined \def \blpage #1{#1}\fi
\ifx \burl  \undefined \def \burl#1{\textsf{#1}}\fi
\ifx \doiurl  \undefined \def \doiurl#1{\url{https://doi.org/#1}}\fi
\ifx \betal  \undefined \def \betal{\textit{et al.}}\fi
\ifx \binstitute  \undefined \def \binstitute#1{#1}\fi
\ifx \binstitutionaled  \undefined \def \binstitutionaled#1{#1}\fi
\ifx \bctitle  \undefined \def \bctitle#1{#1}\fi
\ifx \beditor  \undefined \def \beditor#1{#1}\fi
\ifx \bpublisher  \undefined \def \bpublisher#1{#1}\fi
\ifx \bbtitle  \undefined \def \bbtitle#1{#1}\fi
\ifx \bedition  \undefined \def \bedition#1{#1}\fi
\ifx \bseriesno  \undefined \def \bseriesno#1{#1}\fi
\ifx \blocation  \undefined \def \blocation#1{#1}\fi
\ifx \bsertitle  \undefined \def \bsertitle#1{#1}\fi
\ifx \bsnm \undefined \def \bsnm#1{#1}\fi
\ifx \bsuffix \undefined \def \bsuffix#1{#1}\fi
\ifx \bparticle \undefined \def \bparticle#1{#1}\fi
\ifx \barticle \undefined \def \barticle#1{#1}\fi
\bibcommenthead
\ifx \bconfdate \undefined \def \bconfdate #1{#1}\fi
\ifx \botherref \undefined \def \botherref #1{#1}\fi
\ifx \url \undefined \def \url#1{\textsf{#1}}\fi
\ifx \bchapter \undefined \def \bchapter#1{#1}\fi
\ifx \bbook \undefined \def \bbook#1{#1}\fi
\ifx \bcomment \undefined \def \bcomment#1{#1}\fi
\ifx \oauthor \undefined \def \oauthor#1{#1}\fi
\ifx \citeauthoryear \undefined \def \citeauthoryear#1{#1}\fi
\ifx \endbibitem  \undefined \def \endbibitem {}\fi
\ifx \bconflocation  \undefined \def \bconflocation#1{#1}\fi
\ifx \arxivurl  \undefined \def \arxivurl#1{\textsf{#1}}\fi
\csname PreBibitemsHook\endcsname

\bibitem[\protect\citeauthoryear{D{\k{a}}browski and Sitarz}{2015}]{DS15}
\begin{barticle}
\bauthor{\bsnm{D{\k{a}}browski}, \binits{L.}},
\bauthor{\bsnm{Sitarz}, \binits{A.}}:
\batitle{An asymmetric noncommutative torus}.
\bjtitle{Symmetry, Integrability and Geometry: Methods and Applications
  (SIGMA)}
\bvolume{11},
\bfpage{075}
(\byear{2015})
\doiurl{10.3842/SIGMA.2015.075}
\end{barticle}
\endbibitem

\bibitem[\protect\citeauthoryear{D\k{a}browski et~al.}{2023}]{DSZ23}
\begin{barticle}
\bauthor{\bsnm{D\k{a}browski}, \binits{L.}},
\bauthor{\bsnm{Sitarz}, \binits{A.}},
\bauthor{\bsnm{Zalecki}, \binits{P.}}:
\batitle{Spectral metric and Einstein functionals}.
\bjtitle{Advances in Mathematics}
\bvolume{427},
\bfpage{109128}
(\byear{2023})
\doiurl{10.1016/j.aim.2023.109128}
\end{barticle}
\endbibitem

\bibitem[\protect\citeauthoryear{Connes}{1994}]{Co94}
\begin{bbook}
\bauthor{\bsnm{Connes}, \binits{A.}}:
\bbtitle{Noncommutative Geometry},
\bedition{1st} edn.
\bpublisher{Academic Press},
\blocation{San Diego CA U.S.A.}
(\byear{1994}).
\bcomment{ISBN: 978-0121858605}
\end{bbook}
\endbibitem

\bibitem[\protect\citeauthoryear{Connes and Tretkoff}{2011}]{CoTr11}
\begin{botherref}
\oauthor{\bsnm{Connes}, \binits{A.}},
\oauthor{\bsnm{Tretkoff}, \binits{P.}}:
The Gauss-Bonnet theorem for the noncommutative two torus, Noncommutative
  geometry, arithmetic, and related topics, 141-158.
Johns Hopkins Univ. Press, Baltimore, MD
(2011).
\url{https://www.jstor.org/stable/43302854}
\end{botherref}
\endbibitem

\bibitem[\protect\citeauthoryear{Connes and Moscovici}{2014}]{CoMo14}
\begin{barticle}
\bauthor{\bsnm{Connes}, \binits{A.}},
\bauthor{\bsnm{Moscovici}, \binits{H.}}:
\batitle{Modular curvature for noncommutative two-tori}.
\bjtitle{Journal of the American Mathematical Society}
\bvolume{27}(\bissue{3}),
\bfpage{639}--\blpage{684}
(\byear{2014})
\doiurl{10.1090/S0894-0347-2014-00793-1}
\end{barticle}
\endbibitem

\bibitem[\protect\citeauthoryear{D{\k{a}}browski et~al.}{2024}]{DSZ24}
\begin{barticle}
\bauthor{\bsnm{D{\k{a}}browski}, \binits{L.}},
\bauthor{\bsnm{Sitarz}, \binits{A.}},
\bauthor{\bsnm{Zalecki}, \binits{P.}}:
\batitle{{Spectral torsion}}.
\bjtitle{Commun. Math. Phys.}
\bvolume{405}(\bissue{5}),
\bfpage{130}
(\byear{2024})
\doiurl{10.1007/s00220-024-04950-7}
\end{barticle}
\endbibitem

\bibitem[\protect\citeauthoryear{D\k{a}browski et~al.}{2024}]{DSZ24a}
\begin{barticle}
\bauthor{\bsnm{D\k{a}browski}, \binits{L.}},
\bauthor{\bsnm{Sitarz}, \binits{A.}},
\bauthor{\bsnm{Zalecki}, \binits{P.}}:
\batitle{Spectral metric and Einstein functionals for the Hodge - Dirac
  operator}.
\bjtitle{Journal of Noncommutative Geometry}
(\byear{2024})
\doiurl{10.4171/jncg/573}
\end{barticle}
\endbibitem

\bibitem[\protect\citeauthoryear{Wodzicki}{1987}]{Wo87}
\begin{bbook}
\bauthor{\bsnm{Wodzicki}, \binits{M.}}:
In: \beditor{\bsnm{Manin}, \binits{Y.I.}} (ed.)
\bbtitle{Noncommutative residue. {C}hapter I. {F}undamentals},
pp. \bfpage{320}--\blpage{399}.
\bpublisher{Springer},
\blocation{Berlin, Heidelberg}
(\byear{1987}).
\doiurl{10.1007/BFb0078372}
\end{bbook}
\endbibitem

\bibitem[\protect\citeauthoryear{Lesch}{2017}]{Le17}
\begin{barticle}
\bauthor{\bsnm{Lesch}, \binits{M.}}:
\batitle{Divided differences in noncommutative geometry: rearrangement lemma,
  functional calculus and expansional formula}.
\bjtitle{Journal of Noncommutative Geometry}
\bvolume{11}(\bissue{1}),
\bfpage{193}--\blpage{223}
(\byear{2017})
\doiurl{10.4171/JNCG/11-1-6}
\end{barticle}
\endbibitem

\bibitem[\protect\citeauthoryear{Hekkelman et~al.}{2024}]{HeMDvN24}
\begin{botherref}
\oauthor{\bsnm{Hekkelman}, \binits{E.-M.}},
\oauthor{\bsnm{McDonald}, \binits{E.}},
\oauthor{\bsnm{Nuland}, \binits{T.D.}}:
Multiple operator integrals, pseudodifferential calculus, and asymptotic
  expansions.
(2024)
{\href{https://arxiv.org/abs/2404.16338}{{arXiv:2404.16338}}}
\end{botherref}
\endbibitem

\bibitem[\protect\citeauthoryear{Fathizadeh and Khalkhali}{2012}]{FaKh12}
\begin{barticle}
\bauthor{\bsnm{Fathizadeh}, \binits{F.}},
\bauthor{\bsnm{Khalkhali}, \binits{M.}}:
\batitle{The Gauss--Bonnet theorem for noncommutative two tori with a general
  conformal structure}.
\bjtitle{Journal of noncommutative geometry}
\bvolume{6}(\bissue{3}),
\bfpage{457}--\blpage{480}
(\byear{2012})
\doiurl{10.4171/JNCG/97}
\end{barticle}
\endbibitem

\bibitem[\protect\citeauthoryear{Fathizadeh and Khalkhali}{2013}]{FaKh13}
\begin{barticle}
\bauthor{\bsnm{Fathizadeh}, \binits{F.}},
\bauthor{\bsnm{Khalkhali}, \binits{M.}}:
\batitle{Scalar curvature for the noncommutative two torus}.
\bjtitle{Journal of Noncommutative Geometry}
\bvolume{7}(\bissue{4}),
\bfpage{1145}--\blpage{1183}
(\byear{2013})
\doiurl{10.4171/JNCG/145}
\end{barticle}
\endbibitem

\bibitem[\protect\citeauthoryear{Fathizadeh and Khalkhali}{2011}]{FaKh13a}
\begin{botherref}
\oauthor{\bsnm{Fathizadeh}, \binits{F.}},
\oauthor{\bsnm{Khalkhali}, \binits{M.}}:
Weyl's law and Connes' trace theorem for noncommutative two tori.
Letters in Mathematical Physics
\textbf{103}
(2011)
\doiurl{10.1007/s11005-012-0593-2}
\end{botherref}
\endbibitem

\bibitem[\protect\citeauthoryear{Rosenberg}{2013}]{Ro13}
\begin{barticle}
\bauthor{\bsnm{Rosenberg}, \binits{J.}}:
\batitle{Levi-Civita's theorem for noncommutative tori}.
\bjtitle{SIGMA. Symmetry, Integrability and Geometry: Methods and Applications}
\bvolume{9},
\bfpage{071}
(\byear{2013})
\doiurl{10.3842/SIGMA.2013.071}
\end{barticle}
\endbibitem

\bibitem[\protect\citeauthoryear{Floricel et~al.}{2019}]{FGK19}
\begin{barticle}
\bauthor{\bsnm{Floricel}, \binits{R.}},
\bauthor{\bsnm{Ghorbanpour}, \binits{A.}},
\bauthor{\bsnm{Khalkhali}, \binits{M.}}:
\batitle{The Ricci curvature in noncommutative geometry.}
\bjtitle{Journal of Noncommutative Geometry}
(\byear{2019})
\doiurl{10.4171/jncg/324}
\end{barticle}
\endbibitem

\end{thebibliography}
\end{document}